\documentclass[10pt]{article}
\usepackage{amssymb,amsbsy,amsmath,amsfonts}
\usepackage{latexsym,euscript,exscale}

\title{\bf A topological version of the Bergman property}
\author {\textsc{Christian Rosendal}}
\date {}
\newcommand {\ca} {{2^\N}}

\newcommand {\F}{\mathbb F}

\newcommand {\AAA}{\mathbb A}
\newcommand {\N}{\mathbb N}
\newcommand {\Q}{\mathbb Q}
\newcommand {\R}{\mathbb R}
\newcommand {\Z}{\mathbb Z}
\newcommand {\C}{\mathbb C}
\newcommand {\U}{\mathbb U}

\newcommand {\HH}{\mathbb H}

\newcommand{\q}{\fed Q}
\newcommand{\V}{\mathbb V}
\newcommand{\W}{\mathbb W}

\newcommand{\K}{\mathsf K}

\newcommand{\norm}[1]{\lVert#1\rVert}

\newcommand{\normal}{\trianglelefteq}
\newcommand{\om}{\omega}
\newcommand{\eps}{\epsilon}

\newcommand{\sym}{\vartriangle}
\newcommand{\tom} {\emptyset}

\newcommand{\begr}{\!\upharpoonright}

\newcommand{\saa}{\Longrightarrow}

\newcommand{\til}{\rightarrow}

\newcommand{\Lim}[1]{\mathop{\longrightarrow}\limits_{#1}}

\newcommand {\del}{ \; \big| \;}

\newcommand {\for}{\bigcup}

\newcommand {\go} {\mathfrak}
\newcommand {\ku} {\mathcal}

\newcommand{\ov}{\overline}
\newcommand{\inv}{^{-1}}

\newcommand {\e} {\exists}
\renewcommand {\a} {\forall}

\newcommand{\fed}{\boldsymbol}

\newcommand{\pf}{

\smallskip

\noindent {\it Proof : }}

\newcommand{\pff}{$\hfill  \Box$

\smallskip }

\newtheorem{thm}{Theorem}[section]
\newtheorem{cor}[thm]{Corollary}
\newtheorem{lemme}[thm]{Lemma}
\newtheorem{prop} [thm] {Proposition}
\newtheorem{defi} [thm] {Definition}

\newtheorem{quest}[thm]{Question}

\newtheorem{claim}[thm] {Claim}

\newtheorem{nota}[thm]{Notation}

\begin{document}
\maketitle
\begin{abstract}
A topological group $G$ is defined to have property (OB) if any $G$-action by isometries on a metric space, which is separately continuous, has bounded orbits. We study this topological analogue of the socalled Bergman property in the context of Polish groups, where we show it to have several interesting reformulations and consequences. We subsequently apply the results obtained in order to verify property (OB) for a number of groups of isometries and homeomorphism groups of compact metric spaces. We also give a proof that the isometry group of the rational Urysohn metric space of diameter $1$ is Bergman.\footnotemark
\end{abstract}

\section{Introduction}
\footnotetext{MSC primary 03E15, secondary 22F05}
We study in this paper a topological version of the Bergman property. This latter property has been the object of intense
scrutiny by a number of people, since it was first discovered to hold for the infinite symmetric group, $\ku S_\infty$, by George
M. Bergman \cite{ber} sometime in 2003.

\begin{defi} A group $G$ is said to have the {\em Bergman property} (or is {\em Bergman}) if whenever $W_0\subseteq
W_1\subseteq \ldots\subseteq G=\for_n W_n$, there are $n$ and $k$ such that $G=W_n^k$.\end{defi}

We now have a large number of interesting results concerning this property, but most surprising is perhaps the fact that many large
permutation groups are indeed Bergman. One can pick some of these results from \cite{corn,drogob,drohol,kecros,mil,tol1,tol2}. For
example, it holds for automorphism groups of $2$-transitive linear orders (Droste and Holland \cite{drohol}), the group of measure
preserving automorphisms of the unit interval (Miller \cite{mil}) and oligomorphic permutation groups with ample generics, e.g.,
automorphism groups of $\om$-stable, $\om$-categorical structures (Kechris and Rosendal \cite{kecros}).

To see what the Bergman property is really worth, it is useful to consider some of its consequences and reformulations. First of all, it is clear that no Bergman group can be written as a union of a countable chain of proper subgroups, or, in other words, Bergman groups have {\em uncountable cofinality}. Similarly, if $1\in E=E\inv$ is a generating set for a Bergman group, then there is some finite power $n$ such that every element of the group can be written as a word of length $n$ in $E$. We express this by saying that the group is {\em Cayley bounded}, since it corresponds to every Cayley graph being of bounded diameter with respect to the word metric. Both uncountable cofinality and Cayley boundedness have been studied in the literature, though aparently mostly independently of eachother. Uncountable cofinality grew out of J.P. Serre's work \cite{serre} on property (FA), which is a fixed point property for actions on trees, in which it proved to be one of the three conditions in his equivalent formulation of property (FA) for uncountable groups. It was first proved to hold for certain profinite groups by  S. Koppelberg and J. Tits \cite{koptits} and has subsequently been verified for a large number of primarily subgroups of the infinite symmetric group. S. Shelah \cite{shelah}, on the other hand, has constructed a group of cardinality $\aleph_1$ having width $240$ with respect to any generating set, so, in particular, the group is Cayley bounded, but moreover, has no uncountable proper subgroups. Thus, Shelah's example also has uncountable cofinality. As a matter of fact, as was noticed by Droste and Holland \cite{drohol}, these two properties together are equivalent to being Bergman, so Shelah's example is Bergman too.

However, perhaps more useful for the geometric theory of Bergman groups is the following basic characterisation, of which I learned the equivalence of (2) and (3) from B.D. Miller \cite{mil} and where the equivalence of (1) and (3) was independently noticed by Y. de Cornulier
\cite{corn} and V. Pestov.
\begin{thm} The following conditions are equivalent for a group $G$.\\
(1) Whenever $G$ acts by isometries on a metric space $(X,d)$ every orbit is bounded. \\
(2) Any left-invariant metric on $G$ is bounded. \\
(3) $G$ is Bergman. \\
(4) Whenever $G$ acts on a metric space $(X,d)$ by mappings, which are Lipschitz for large distances, every orbit is bounded. \\
(5) Whenever $G$ acts by uniform homeomorphisms on a geodesic space $(X,d)$ every orbit is bounded.
\end{thm}

Now, of course, (1),(4) and (5) are really  properties one would tend to study in connection with topological groups modulo some continuity condition, and these are indeed the main objects of the present paper. We therefore propose the following definition.

\begin{defi}A topological group $G$ is said to have {\em property (OB)} if whenever $G$ acts by isometries on a metric space $(X,d)$, such that for every $x\in X$ the function $g\in G\mapsto gx\in X$ is continuous, then every orbit is bounded.\end{defi}

Similar properties have previously been considered by Jan Hejcman \cite{hej1} in 1959 (see also his recent paper \cite{hej2}) and
later by Christopher Atkin \cite{atkin} under the name of {\em boundedness} in the context of uniform spaces. However, let me point of the differences between their notions and property (OB). A topological group $G$ is {\em bounded} if for any non-empty open subset $U\subseteq G$ there is a finite set $A\subseteq G$ and a number $n$ such that $G=U^nA$. This, however, is easily seen to fail for our Bergman group par excellence, $\ku S_\infty$. For if we choose $U$ to be an open sugroup of denumerable index, as for example the isotropy subgroup of $0\in \N$, then clearly $U^nA=UA\neq \ku S_\infty$ for all $n$ and finite $A$. Nevertheless, the two notions do turn out to be equivalent in the context of abelian groups, but most of the groups considered here are very much non-abelian.

One of the reasons for our interest in the property comes from the fact that it can be seen as an addition to
a well-known spectrum of properties studied in geometric group theory, namely properties (FA), (FH), (T), amenability, etc. One
easily sees that property (OB) implies property (FH) and actually, we shall see that it provides a fairly comprehensive class of new
examples of non-locally compact groups with property (FH).

Much of the work on these properties has been restricted to the locally compact setting, where the strongest tools are available
(e.g., Haar measure). But over the years, a number of very interesting results concerning the dynamics of non-locally compact
Polish groups have surfaced, for example on the unitary group of $\ell_2$, where  Gromov and Milman proved that it is extremely
amenable \cite{gromil} and Bekka proved that it has property (T) \cite{bek}. Moreover, in logic, where, e.g., automorphism groups
of countable structures tend to be non-locally compact, there is a multitude of results on permutation groups, e.g., Truss
\cite{truss} and Hodges, Hodkinson, Lascar, and Shelah \cite{hoholash}, and also for more inclusive classes of Polish
groups, e.g., Becker and Kechris \cite{beckec} and Hjorth \cite{greg}.

Thus for Polish groups it is natural to look for a similar characterisation of property (OB), and indeed we have the following result.
\begin{thm}\label{characterisation}The following are equivalent for a Polish group $G$:\\
(i) $G$ has property (OB).\\
(ii) Whenever $W_0\subseteq W_1\subseteq \ldots\subseteq G$ is an increasing exhaustive sequence of sets with the Baire property, there are $n$ and $k$ such that $G=W_n^k$.\\
(iii) Any compatible left-invariant metric on $G$ is bounded.\\
(iv) $G$ is finitely generated of bounded width over any non-empty open subset.
\end{thm}
For example, a locally compact Polish group has property (OB) if and only if it is compact.

This gives an indication of how to think of these properties. Namely, one should think of the Bergman property as a strong
generalisation of finiteness and of property (OB) as a strong generalisation of compactness. Surprisingly though, this
"compactness" can, apart from compact groups, only be found in large, very much non-locally compact groups.

We study first the dynamics of property (OB) groups acting continuously by H\"older mappings, showing that in this case the
closure of the orbits gives a decomposition of the phase space into pieces on which the group acts minimally (often
denoted by {\em semi-simplicity}). And secondly we consider the closure properties of the class of property (OB)
groups, for example, it is quite easily seen that it is closed under infinite products, group extensions over a property (OB) group, and behaves well with respoect to short exact sequences. I.e., if $\pi: G\til H$ is a continuous homomorphism with dense image where $G$ has property (OB), then so does $H$. Most interesting in this connection is the fact that it passes to subgroups of finite index, for which we give a geometric proof. In his paper \cite{ber}, Bergman originally asked whether his property was preserved between a group and its subgroups of finite index, and A. Kh\'elif, in an announcement \cite{khe}, stated that this is indeed the case. However, the mentioned geometric proof, which works also for the Bergman property, shows the usefulness of the reformulation of the Bergman property in terms of isometric actions, where one has the added advantage of geometric intuition.

The known examples of Bergman groups have mostly been groups of symmetries of various countable structures of spaces, and, maybe apart from the obvious counter-examples of profinite groups of which examples can be found in Saxl, Shelah and Thomas \cite{saxl}, Thomas \cite{thomas} and de Cornulier \cite{corn}, one has the feeling that this is the most likely place to discover these. Somehow there is an intuition that Bergman groups are not constructed from below beginning with the simplest groups and using algebraic constructions such as direct sums and group extensions. This feeling seems to be reinforced by the following result.
\begin{thm}
A solvable group is Bergman if and only if it is finite.
\end{thm}

However, a large chunk of the paper is concerned with the verification and construction of groups which are either Bergman
or have property (OB). We consider a fair number of examples, beginning with groups connected with the unit circle. We subsequently turn Theorem  \ref{characterisation} on its head and instead ask for when the isometry group of a bounded complete metric space has property (OB). We provide one suffificient condition that is also of independent interest and use this to show that the isometry group of the Urysohn metric space of diameter $1$ has property (OB). However, in the case of the rational Urysohn metric space of diameter $1$, one can take advantage of the recent deep results of S. Solecki from \cite{sol}, and use this to show outright that its isometry group is Bergman. 

\begin{thm} Let $\U_1$ be the Urysohn metric space of diameter $1$ and $\Omega$ the rational Urysohn metric space of diameter $1$. Then ${\rm Iso}(\U_1)$ has property (OB) and ${\rm Iso}(\Omega)$ the Bergman property.\end{thm}

We then study a model theoretic version of the unitary group of $\ell_2$ in some depth. This is a subgroup $\ku U(\V)$ that sits as a dense subgroup in $\ku U(\ell_2)$, which we prove to have ample generics, the main tool used in proving the Bergman property for automorphism groups of $\om$-stable, $\om$-categorical structures in \cite{kecros} and previously introduced by Hodges, Hodkinson, Lascar, and Shelah \cite{hoholash}. From ample generics, we prove that also $\ku U(\V)$ is Bergman and has a number of other properties, e.g., the small index property and satisfies automatic continuity of homomorphisms. We use the analysis of $\ku U(\V)$ to further study the dynamics of actions by H\"older mappings $\ku U(\ell_2)$, obtaining the strange conclusion that these are to some extent determined by the action of any bilateral shift on $\ell_2$.

Our final collection of examples comes from topology, where we prove property (OB) for homeomorphism groups of spheres and of the Hilbert cube.

\begin{thm} ${\rm Hom}(S^m)$ and ${\rm Hom}(Q)$ have property (OB) and, by consequence, any Polish group is a subgroup of a Polish group with property (OB).\end{thm}

In the final section of the paper we consider property (FA), and provide a simple proof of a result of Dugald Macpherson and Simon Thomas
stating that if a Polish group with a comeagre conjugacy class acts on a tree, there every element of the group fixes a vertex or
an edge. Actually, we extend their theorem to all actions on $\Lambda$-trees, though of course, the case of simplicial trees it
ultimately the most interesting, due to the structure theory of Serre for groups with property (FA) \cite{serre}.

\

Though we shall from time to time use a little bit of descriptive set theory, the article should be comprehensible to the general
analyst. Really, all one needs to know is the definition of a {\em Polish space} (a completely metrisable separable topological
space), a {\em Polish group} (a topological group whose topology is Polish), {\em Borel sets} (the sets belonging to the $\sigma$-algebra generated by the open sets), {\em analytic set} (a subset in a Polish space which is the continuous image of a Borel set) and sets having the {\em Baire property} (i.e. sets $A$ such that for some open set $U$ and some meagre set $M$, $A=U\sym M$). A basic result of Lusin and Sierpi\'nski says that analytic sets have the Baire property.

\

I am grateful to Alekos Kechris for many welcome criticisms, Benjamin Miller for first getting me interested in Bergman groups
and to Pandelis Dodos for suggesting that abelian groups shouldn't be Bergman. Vladimir Pestov, who personally has also been
interested in the subject of this paper, helped me with different pointers to the literature. Finally, I thank Stevo Todor\v cevi\'c
for a late night discussion at Caltech.

\section{The Bergman property}
\begin{defi} A group $G$ is said to have the {\em Bergman
property} if whenever $W_0\subseteq W_1\subseteq \ldots\subseteq
G=\for_n W_n$, there are $n$ and $k$ such that
$G=W_n^k$.\end{defi}

For the following, recall that a {\em geodesic space} is a metric space such that between any two points $x$ and $y$ there is a path
of length $d(x,y)$. For example, Banach spaces and $\R$-trees are geodesic spaces. Recall also that a mapping $\phi:X\til X$, where $X$ is a metric space, is called {\em Lipschitz for large distances } if there are constants $c, K$ such that for all $x,y\in X$, $d(\phi x,\phi y)\leq c\cdot d(x,y)+K$.

The following result states the basic equivalent formulations of the Bergman property. The equivalence of (2) and (3) I learned from B.D. Miller
\cite{mil} and the equivalence of (1) and (3) was independently noticed by Y. de Cornulier \cite{corn} and V. Pestov.
\begin{thm}\label{basic} The following conditions are equivalent for a group $G$.\\
(1) Whenever $G$ acts by isometries on a metric space $(X,d)$ every orbit is bounded. \\
(2) Any left-invariant metric on $G$ is bounded. \\
(3) $G$ is Bergman. \\
(4) Whenever $G$ acts on a metric space $(X,d)$ by mappings which are Lipschitz for large distances, every orbit is bounded. \\
(5) Whenever $G$ acts by uniform homeomorphisms on a geodesic space $(X,d)$ every orbit is bounded.
\end{thm}

\pf Clearly, $1\saa 2$ is trivial.

$2\saa 3$: Suppose that $W_0\subseteq W_1\subseteq W_2\subseteq\ldots \subseteq G$ is an exhaustive
sequence of subsets of $G$. Notice that then $W_0\cap W_0\inv \subseteq W_1\cap W_1\inv \subseteq \ldots \subseteq G$ is also
exhaustive. So we can suppose that the $W_n$ are symmetric, and by renumbering the sequence, we can also suppose that $W_0=\{1\}$.
Notice now that the following left-invariant metric on $G$ is bounded if and only if $G=W_n^k$ for some $n$ and $k$:
$$
d(f,g)=\min (k_1+k_2+\ldots+k_n\del \e h_i\in W_{k_i}\; fh_1\ldots h_n=g).
$$

$3\saa 4$: Assume now that 3 holds and that $G$ acts on a metric space $(X,d)$ by mappings which are Lipschitz for large distances and find for each $g\in G$ constants $c_g$ and $K_g$ witnessing this. Then, if
$g_1,\ldots,g_k\in G$ and $c_{g_1},\ldots, c_{g_k},K_{g_1},\ldots, K_{g_k}\leq M$, 
\begin{equation}\begin{split}
d(g_1\ldots g_k\cdot x,g_1\ldots g_k\cdot y)&\leq M\cdot
d(g_2\ldots g_k\cdot x,g_2\ldots g_k\cdot
y)+M\\
&\leq M^2\cdot d(g_3\ldots g_k\cdot
x,g_3\ldots g_k\cdot y)+M^2+M\\
&\leq\ldots\\
&\leq  M^k \cdot d(x,y) +M^k+M^{k-1}+\ldots+M\\
&\leq M^k\cdot (d(x,y)+k)
\end{split}\end{equation}
Now, fix an $x_0\in X$ and let for $n\geq 1$
$$
W_n=\{g\in G\del c_g,K_g\leq n\;\&\; d(x_0,g\cdot x_0)\leq n\}
$$
This is clearly an increasing exhaustive sequence of subsets of
$G$, so for some $M$ and $k$, $G=W_M^k$. We claim that $x_0$'s orbit is bounded in diameter by $2k^2M^k$. For if $g=g_1\ldots g_k\in G$, with $g_1,\ldots, g_k\in W_M$, 
\begin{equation}\begin{split}
d(x_0,g\cdot x_0)&=d(x_0,g_1\ldots g_k\cdot x_0)\\
&\leq d(x_0,g_1\cdot x_0)+d(g_1\cdot x_0, g_1g_2\cdot
x_0)+\ldots+d(g_1\ldots g_{k-1}\cdot x_0,g_1\ldots g_{k-1}g_k\cdot
x_0)\\
&\leq d(x_0,g_1\cdot x_0)+ M (d(x_0, g_2\cdot
x_0)+1)+M^2(d(x_0,g_3\cdot x_0)+2)+\\
&\quad\;\ldots+M^{k-1}(d(x_0,g_k\cdot
x_0)+k-1)\\
&\leq M+M(M+1)+M^2(M+2)+\ldots+M^{k-1}(M+k-1)\\
&\leq 2k^2M^k
\end{split}\end{equation}
So if $x$ is any other point of $X$, then
\begin{equation}\begin{split}
d(g\cdot x, x)&\leq d(g\cdot x, g\cdot x_0)+d(g\cdot
x_0,x_0)+d(x_0,x)\\
&\leq M^k(d(x,x_0)+k)+2k^2M^k+d(x,x_0)
\end{split}\end{equation}
Whence $x$'s orbit is bounded and thus showing 4.

$4\saa 5$: This implication follows from the general fact that any uniformly continuous mapping $\phi$ on a geodesic space $(X,d)$ is Lipschitz for large distances. To see this, notice that for some $\eps>0$ and all $x,y\in X$,
$$
d(x,y)\leq\eps\til d(\phi x,\phi y)\leq1.
$$
So if $d(x,y)\leq N\cdot \eps$, there are, since $(X,d)$ is
geodesic, $x_0=x,x_1,\ldots,x_N=y\in X$ with $d(x_i,x_{i+1})\leq
\eps$, whence
$$
d(\phi x,\phi y)\leq d(\phi x_0,\phi x_1)+d(\phi
x_1,\phi x_2)+\ldots+d(\phi x_{N-1},\phi x_N)\leq N.
$$
In other words, for all $x,y\in X$,
$$
d(\phi x,\phi y)\leq \Big\lceil\frac{d(x,y)}{\eps}\Big\rceil<\frac1\eps d(x,y)+1.
$$

$5\saa 1$: Suppose that $G$ acts by isometries on a metric space $(X,d)$. Now, $(X,d)$ might not be geodesic, but can be extended to a geodesic space as follows:

For any $x,y\in X$, let $\alpha(x,y)$ be a distinct isometric copy of $[0,d(x,y)]$. We let $\tilde X$ be the quotient space of $\mathsf X=\for_{x,y\in X}\alpha(x,y)$ obtained by for all $x,y,z\in X$ identifying the left endpoints of $\alpha(x,y)$ and $\alpha(x,z)$, identifying the right endpoints of $\alpha(y,x)$ and $\alpha(z,x)$, and identifying the right endpoint of $\alpha(y,x)$ with the left endpoint of $\alpha(x,z)$. 
We define a metric $\tilde d$ on $\tilde X$ as follows: For $a \in\alpha(x,y)$ and $b\in\alpha(z,u)$, put
\begin{displaymath}\tilde d(a,b)=\min\left\{
\begin{array}{l}
|a-0|+d(x,z)+|0-b|\\
|a-0|+d(x,u)+|d(z,u)-b|\\
|a-d(x,y)|+d(y,z)+|0-b|\\
|a-d(y,u)|+d(y,u)+|d(z,u)-b|
\end{array}\right.
\end{displaymath}
Then $(\tilde X, \tilde d)$ contains $(X,d)$ isometrically (sending $x\in X$ to the equivalence class of the unique (end)point of $\alpha(x,x)$). Moreover, the isometric action of $G$ on $X$ extends to $\tilde X$ by letting $g\in G$ send $\alpha(x,y)$ isometrically and orderpreservingly to $\alpha(gx,gy)$. Thus, $G$ acts by isometries on the geodesic space $\tilde X$ and hence, if $5$ holds, then every orbit of $\tilde X$ and hence every orbit of $X$ is bounded.
\pff

This allows us to give the following nice proof that the Bergman property is preserved under short exact sequences: Clearly, if
$H\normal G$ and any action by isometries of $G$ has bounded orbits, then any action by isometries of $G/H$ has bounded orbits.
Conversely, assume any action by isometries of $G/H$ and of $H$ has bounded orbits and that $G$ acts by isometries on $(X,d)$. Let
$\ku O$ be the closure of an $H$ orbit in $X$ and let $A=\{g\cdot \ku O\del g\in G\}$ be equipped with the Hausdorff metric $d_H$.
Then $G$ acts transitively by isometries on $(A,d_H)$ and the action factors through $G/H$. Therefore, $A$ is bounded and so any
orbit is bounded in $X$.

\section{Polish groups with property (OB)}
\begin{defi}A topological group $G$ is said to have {\em property
(OB)} if whenever $G$ acts by isometries on a metric space
$(X,d)$, such that for every $x\in X$ the mapping $g\in G\mapsto
gx\in X$ is continuous, then every orbit is bounded.\end{defi}

The preceding section is exclusively concerned with discrete groups, but we shall see that in the case of Polish groups there
are again nice equivalent formulations of property (OB). We first notice that property (OB) can be slightly reformulated
for Polish groups.

\begin{lemme}\label{joint}Let $G$ be a Polish group acting by
homeomorphisms on a metrisable space $X$, such that the mapping
$g\in G\mapsto gx\in X$ is continuous for every $x\in X$. Then the
action is actually jointly continuous, i.e., $(g,x)\mapsto g\cdot
x$ is continuous from $G\times X$ to $X$.\end{lemme}

\pf Assume that $g_n\til g$ and $x_n\til x$. Since the mapping
$h\in G\mapsto hy\in X$ is continuous for every $y\in X$, we see
that $G\cdot y$ is a continuous image of a separable space and
thus separable for every $y$. Hence $Y=G\cdot x\cup \for_n G\cdot
x_n$ is an invariant separable subspace of $X$. Moreover, the
action of $G$ on $Y$ is separately continuous, so, as $Y$ is
metrisable, the action of $G$ on $Y$ is jointly continuous
(Kechris \cite{kec} (9.16)). Therefore, $g_nx_n\til gx$
in $Y$ and thus also in $X$.\pff

In particular, a Polish group has property (OB) if and only if all
of its continuous actions by isometries on separable metric spaces have
bounded orbits.

Property (OB) and the Bergman property fit quite nicely into the well-known hierarchy of group theoretical fixed-point
properties such as property (T), (FH), (FA) etc. As first sight they do not appear to be a fixed-point properties, but it all depends on the perspective, as, for example, the Bergman property is equivalent to a fixed point property for its induced actions on the hyperspace of bounded subsets of any metric space it acts upon. 

\begin{defi} A group $G$ is said to have {\em property
(FA)} if whenever it acts by automorphisms on a combinatorial tree
(i.e. a uniquely path connected graph) it either fixes a
vertex or an edge.

A topological group $G$ has {\em property (topFA)} if whenever it
acts by automorphisms on a combinatorial tree, such that the
stabilisers of vertices are open, then it fixes either a vertex or
an edge.

A group $G$ is said to have {\em property (algFH)} if whenever it
acts by isometries on a real Hilbert space $\ku H$, then it fixes
a vector.

A topological group $G$ is said to have {\em property (FH)} if
whenever it acts by isometries on a real Hilbert space $\ku H$
such for all $\xi\in \ku H$ the mapping $g\in G\mapsto g\cdot
\xi\in \ku H$ is continuous, then it fixes a vector.
\end{defi}

Admittedly, the fixed point property on trees is mainly interesting in its algebraic version, property (FA). Indeed, it is the main object of Serre's book \cite{serre} in which he shows that it is equivalent to the conjunction of (i) the group has no infinite cyclic quotients, (ii) the group is not a non-trivial free product with amalgamation and (iii) the group is not the union of a countable chain of proper subgroups.  The fixed point property on Hilbert spaces has correspondingly mostly been studied for countable discrete groups (in which case properties (algFH) and (FH) coincide) and for locally compact groups, where one is interested in property (FH). It is also well-known that (algFH) is stronger than (FA) and similarly, (FH) is stronger than (topFA). Moreover, for a group of isometries of Hilbert space to fix a point it is enough that there should be a bounded orbit. This follows from the lemma of the centre (see Bekka, de la Harpe and Valette \cite{behava}).
So the following proposition sums up the connections between our properties.

\begin{prop}The following diagram of implications holds for
topological and abstract groups.
\begin{displaymath}
\begin{array}{ccc}
\textrm {Bergman} &\saa& \textrm {Property } (OB)\\
\Downarrow && \Downarrow \\
\textrm {Property } (algFH) &\saa& \textrm {Property } (FH)\\
\Downarrow && \Downarrow \\
\textrm {Property } (FA)&\saa& \textrm {Property } (topFA)
\end{array}\end{displaymath}
\end{prop}

Now in turn, we will show the basic equivalences of the different formulations of property (OB) for Polish groups. The following extracts the basic properties of the usual proof of the Birkhoff-Kakutani metrisation theorem, see, e.g., Hjorth \cite{greg}, Theorem 7.2.

\begin{lemme}\label{birkhoff} Let $G$ be a topological group and $(V_n)_{n\in\Z}$
a neighbourhood basis at the identity consisting of open sets
such that\\
(I) $V_n=V_n\inv$\\
(II) $G=\for_{n\in \Z}V_n$\\
(III) $V_n^3\subseteq V_{n+1}$

Let $\delta(g_1,g_2)=\inf(2^n\del g_1\inv g_2\in V_n)$ and put
$$
d(g_1,g_2)=\inf (\sum_{i=0}^k\delta(h_i,h_{i+1})\del h_0=g_1,
h_k=g_1)
$$
Then
\begin{equation}\label{bika}
\delta(g_1,g_2)\leq 2d(g_1,g_2)\leq 2\delta(g_1,g_2)
\end{equation}
and $d$ is a left-invariant compatible metric on $G$.
\end{lemme}

\begin{defi} Let $G$ be a Polish group. We say that $G$ is {\em topologically Bergman} if whenever
$$
B_0\subseteq B_1\subseteq B_2\subseteq \ldots \subseteq G
$$
is an exhaustive sequence of subsets with the Baire property, then $G=B_n^k$ for some $n$ and $k$. If there is a $k$ which works for
all sequences $(B_n)$, then we say that $G$ is topologically $k$-Bergman.\end{defi}

\begin{thm}The following are equivalent for a Polish group $G$.\\
(i) $G$ has property (OB).\\
(ii) $G$ is topologically Bergman.\\
(iii) Any compatible left-invariant metric on $G$ is bounded.\\
(iv) $G$ is finitely generated of bounded width over any non-empty open subset.
\end{thm}

\pf (iii)$\saa$(ii): Assume that $G$ is not topologically Bergman as witnessed by some exhaustive sequence of subsets with the Baire
property
$$
B_0\subseteq B_1\subseteq B_2\subseteq \ldots \subseteq G.
$$
By considering $B_0\cap B_0\inv\subseteq B_1\cap B_1\inv\subseteq\ldots$ we can assume that the $B_n$ are symmetric. Then, as $G$ is Polish, some $B_n$ must be non-meagre and contain $1_G$, whence $V=int(B_n^2)\neq \tom$ by Pettis' Lemma. Thus
$$
VB_n\subseteq VB_{n+1}\subseteq \ldots \subseteq G
$$
is an exhaustive sequence of open sets and $(VB_m)^k\subseteq (B_m^3)^k\neq G$ for all $m\geq n$. Put now $V_m=(VB_{n+m})^{3^m}$
and notice that $(V_m)_{m\in \N}$ is an increasing and exhaustive sequence of open neighbourhoods of the identity satisfying
$V_m^3\subseteq V_{m+1}$. Supplementing this sequence with suitable $V_m$ for $m<0$ we get a neighbourhood basis $(V_m)_{m\in
\Z}$ satisfying the conditions of Lemma \ref{birkhoff}. Moreover, as $V_m\neq G$ for all $G$, the resulting metric $d$ is
left-invariant, compatible, but unbounded.

The proof of the implication (ii)$\saa$(i) can be done as in the proof of $3\saa 4$ in Theorem \ref{basic}, and that (i) implies (iii) is
trivial.

(ii)$\saa$(iv): If $V\subseteq G$ is non-empty open and $\{g_n\}_{n\in \N}$ is dense in $G$, then the sequence
$B_n=g_0V\cup\ldots \cup g_nV$ is increasing and exhaustive. So, if $G$ is topologically Bergman, then $G=B_n^k$ for some $n$ and $k$. But then
$G=\big( V \cup\{g_0,\ldots,g_n\}\big)^{2k}$, showing that $G$ is finitely generated of bounded width over $V$.

(iv)$\saa$(ii): Suppose $G$ is finitely generated of bounded width over any non-empty open set and $B_0\subseteq B_1\subseteq
\ldots\subseteq G$ is an increasing exhaustive sequence of sets with the Baire property. Then some $B_n$ is non-meagre and $V=int
B_n^2\neq \tom$. Find some $g_0,\ldots, g_m\in G$ and $k$ such that $G=\big( V \cup\{g_0,\ldots,g_m\}\big)^{k}$. Then $G=\big(
B_n^2 \cup\{g_0,\ldots,g_m\}\big)^{k}\subseteq B_l^{2k}$ for $l\geq n$ large enough such that $1,g_0,\ldots,g_m\in B_l$. \pff

It was shown in Droste and Holland \cite{drohol} that a group $G$ has the Bergman property if and only if $G$ satisfies the
conjunction of the following two properties: 
\begin{itemize}
\item (Uncountable cofinality) Whenever $H_0\leq H_1\leq \ldots\leq G=\for_nH_n$, then $G=H_n$ for some $n$.
\item (Cayley boundedness) Whenever $1\in E=E\inv$ generates $G$ then $G=E^n$ for some $n$.
\end{itemize}
In the same manner, we can define these concepts for Polish groups.

\begin{defi} Let $G$ be a Polish group. We say that $G$ has {\em uncountable topological cofinality} if $G$ is not the union
of a chain of proper open subgroups (or equivalently, a countable chain of subgroups with the Baire property). $G$ is {\em topologically Cayley bounded} if it has finite width with respect to any analytic generating set.
\end{defi}

\begin{prop} \label{cofinality} A Polish group $G$ has property (OB) if and only if it has uncountable topological cofinality and
is topologically Cayley bounded.\end{prop}

\pf Suppose $G$ has property (OB), $H_0\leq H_1\leq\ldots,\leq G=\for_nH_n$ are open subgroups and $1\in E=E\inv$ is an analytic
set generating $G$. Then by property (OB) applied to the sequences $H_0\subseteq H_1\subseteq \ldots \subseteq G$ and $E\subseteq
E^2\subseteq \ldots \subseteq G$, we see that $G=H_n^k=H_n=(E^n)^k=E^{nk}$ for some $n$ and $k$. Thus, $G$ has
both uncountable topological cofinality and is topologically Cayley bounded.

Conversely, suppose $G$ has uncountable topological cofinality, is topologically Cayley bounded and $W_0\subseteq W_1\subseteq \ldots
\subseteq G=\for_nW_n$ are sets with the Baire property. By considering instead a tail subsequence of the exhaustive sequence
$$
W_0\cap W_0\inv\subseteq W_1\cap W_1\inv \subseteq\ldots\subseteq
G
$$
we can suppose each $W_n$ is non-meagre, symmetric and contains $1$. Thus the sequence $\langle W_0\rangle\leq \langle
W_1\rangle\leq \ldots\leq G$ consists of open subgroups and hence one of the $W_n$ generates $G$. As $W_n$ is symmetric and
non-meagre, $int W_n^2$ is symmetric and non-empty, so $W_n\cdot int W_n^2 \cdot W_n\subseteq W_n^4$ is a symmetric generating open
subset of $G$ containing $1$. So $G=W_n^k$ for some $k$.\pff

In the case of locally compact groups, uncountable topological cofinality is clearly equivalent to compact generation. Moreover, if a locally compact, compactly generated group is also topologically Cayley bounded, then its compact generating set generates by a finite power and hence the group is compact. Conversely, compact groups trivially have property (OB). So property (OB) for locally compact Polish groups is just equivalent with compactness, just as the Bergman property for countable groups is equivalent with finiteness. However, we can actually provide a bit more information.

\begin{prop}
A compact Polish group is topologically $2$-Bergman. 
\end{prop}

\pf Assume that $G$ is compact and that
$$
B_0\subseteq B_1\subseteq B_2\subseteq \ldots \subseteq G
$$
is an exhaustive sequence of subsets with the Baire property. Then there is some $B_{n_0}$ which is non-meagre and hence comeagre in
some open set $Vf$, where $V$ is an open neighbourhood of the identity. Pick some symmetric open set $U\subseteq V$ such that
$U^2\subseteq V$ and $g_1,\ldots, g_m\in G$ such that $G=g_1U\cup\ldots\cup g_mU$. Then if $h_i\in g_iU$, we have
$g_iU=h_i(h_i\inv g_i)U\subseteq h_iU^2\subseteq h_iV$ and thus $G= h_1U^2\cup\ldots\cup h_mU^2= h_1V\cup\ldots\cup h_mV$.
Considering now the sequences $(B_j\cap g_iU)_j$ for each $i=1,\ldots, m$, we find $n_1\geq n_0$ such that $B_{n_1}$ is
non-meagre in each of the $g_iU$. Hence, we can find open sets $W_i\subseteq g_iU$ such that $B_{n_1}$ is comeagre in $W_i$ for
each $i$. Pick now $h_i\in W_i\cap B_{n_1}\subseteq g_iU$ and notice that $G=Gf=h_1Vf\cup\ldots\cup h_mVf$. But as $B_{n_1}$ is
comeagre in both $W_i$ and $Vf$, we have  by Pettis' Lemma $h_iVf\subseteq W_iVf\subseteq B_{n_1}^2$. Thus, $G=B_{n_1}^2$, showing that $G$ is topologically $2$-Bergman. \pff

As a locally compact, non-compact Polish group cannot have property (OB), we know that it must have a compatible left-invariant unbounded metric. But actually we can see that this metric can be chosen to be {\em Heine-Borel}, i.e., such that any bounded closed set is compact.

\begin{prop}(Folklore) Let $G$ be a locally compact Polish group. Then $G$ admits a left-invariant compatible Heine-Borel metric. 
\end{prop}

\pf We start by fixing an open neighbourhood basis at the identity $(U_n)_{n\in \N}$ such that $U_{n+1}^3\subseteq U_n$,
$U_n=U_n\inv$ and $U_0$ being relatively compact. Since $G$ is $\sigma$-compact we can also find an increasing sequence of
symmetric relatively compact sets $(V_n)_{n\in \N}$ such that $V_0=U_0$, $V_n^3\subseteq V_{n+1}$ and $G=\for_{n\in \N} V_n$.
Letting now $V_{-n}=U_n$, we see that the sequence $(V_n)_{n\in \Z}$ satisfies the conditions of Lemma \ref{bika}. Let now $d$ be
the metric given by the lemma, we claim that $d$ is Heine-Borel. For any $d$-bounded set $A$ is $\delta$-bounded and thus there is
some $n$ such that for any $g,h\in A$, $g\inv h \in V_n$. In particular, $A$ is contained in some translate of $V_n$ and thus
relatively compact. Hence if $A$ is closed it is compact, showing that $d$ is Heine-Borel.\pff

We should mention that locally compact Polish groups have {\em complete} left-invariant metrics and hence every left-invariant
metric is complete, see Becker \cite{becker}, section  3.

\noindent$\fed {\rm Remark:}$ We should mention that there are
examples of compact Polish groups not being Bergman. In fact,
Koppelberg and Tits \cite{koptits} prove that if $F$ is a finite
non-trivial group, then $F^\N$ has uncountable cofinality if and
only if $F$ is perfect. Thus as the Bergman property implies
uncountable cofinality, we have compact profinite groups without
the Bergman property. We shall see later that, in fact, no solvable
infinite group can be Bergman.

We also see that the two properties (FH) and (OB) do not coincide.
For there are plenty of examples of locally compact, non-compact
Polish groups with property (FH), but of course without property
(OB).

\begin{quest} To what extent do homeomorphism groups of compact
metric spaces have property (OB)? What about the isometry groups
of bounded Polish metric spaces? (Some special cases will be
verified in the following.)\end{quest}

\begin{defi} Recall that a mapping $f$ between metric spaces $X,d$ and $Y,\delta$ is called  a H\"older$(\alpha)$ map for some
$\alpha>0$ if there is a constant $c\geq 1$ such that
$$
\delta(f(x),f(y))\leq c \cdot d(x,y)^\alpha
$$
for all $x,y\in X$. H\"older$(1)$ mappings are thus simply Lipschitz mappings.\end{defi}

\begin{prop} Let $G$ be a Bergman group acting by H\"older maps on a metric space $(X,d)$. Then the action of $G$ is
semi-simple, i.e., $\{\ov{G\cdot x}\}_{x\in X}$ partitions $X$ into (bounded) invariant pieces each on which $G$ acts minimally.
Moreover, there is an $N$ such that any $g\in G$ is H\"older$(\alpha)$ with constant $N$ for some $\alpha\in [1/N,N]$.

The same holds for Polish groups with property (OB) acting continuously and by H\"older maps on a Polish metric space $(X,d)$.
\end{prop}

\pf First assume that $G$ is Bergman. For each $g\in G$ let $\alpha_g\geq 0$ and $c_g\geq 1$ be such that $g$ is
H\"older($\alpha_g$) with constant $c_g$. Thus 
\begin{equation}\begin{split}\label{holder}
d(g_1\cdot\ldots\cdot g_nx, g_1\cdot\ldots\cdot g_ny)&\leq
c_{g_1}d(g_2\cdot\ldots\cdot g_nx, g_2\cdot\ldots\cdot
g_ny)^{\alpha_{g_1}}\\
&\leq c_{g_1}c_{g_2}^{\alpha_{g_1}}d(g_3\cdot\ldots\cdot g_nx,
g_3\cdot\ldots\cdot g_ny)^{\alpha_{g_1}\alpha_{g_2}}\\
&\leq \ldots\\
&\leq c_{g_1}c_{g_2}^{\alpha_{g_1}}\ldots c_{g_n}^{
\alpha_{g_1}\ldots \alpha_{g_{n-1}}}
d(x,y)^{\alpha_{g_1}\ldots\alpha_{g_n}}
\end{split}\end{equation}
Put now $W_n=\{g\in G\del c_g\leq n \;\&\; \alpha_g\in [1/n,n]\}$ and notice that the sequence $W_n$ is increasing and exhaustive.
So as $G$ is Bergman, $G=W_n^k$ for some $n$ and $k$. By the inequality \ref{holder}, we see that there is a fixed $N$ such that any $g\in G$ is H\"older$(\alpha)$ with constant $N$ for some $\alpha\in [1/N,N]$. Thus we have
\begin{equation}\label{equicontinuity}
\a \eps>0\;\:\e \delta> 0\;\: \a x,y\;\:\a g\;\:\:
\big(d(x,y)<\delta \til d(gx,gy)<\eps\big)
\end{equation}
So suppose $x,y\in X$ and that $y\in \ov {G\cdot x}$. Then for any
$x'\in G\cdot x$ and $\eps>0$, we can find $\delta>0$ as above and
$x''\in G\cdot x$ with $d(x'',y)<\delta$. Now, if $x'=gx''$, then
$d(x',gy)=d(gx'',gy)<\eps$, so $x'\in \ov {G\cdot y}$, showing
that $\ov {G\cdot x}=\ov {G\cdot y}$. Thus in $\ov {G\cdot x}$
every orbit is dense, and hence if $u\in\ov {G\cdot x}\cap\ov
{G\cdot z}$ for any $u,z$, then $\ov {G\cdot x}=\ov {G\cdot z}=\ov
{G\cdot u}$. Therefore, $\{\ov{G\cdot x}\}_{x\in X}$ partitions
$X$ and $G$ acts minimally on each piece of the partition. The
usual argument, as in the proof of Theorem \ref{basic}, will also
show that every orbit is bounded.

Now we only have to indicate the proof in the case that $G$ is a
Polish group with property (OB). In this case we fix a countable
dense set $\{x_m\}$ in $X$ and define $W_n$ by
$$
W_n=\{g\in G\del \e \alpha \in [1/n,n]\; \a m,l\; d(gx_m,gx_l)\leq
n\cdot d(x_m,x_l)^\alpha\}
$$
Then $W_n$ is analytic (in fact closed) and hence has the Baire
property.  Notice also that if $g\in W_n$ then $g$ is indeed
H\"older($\alpha$) with constant $n$ for some $\alpha\in [1/n,n]$.
We can thus proceed as before using that $G$ is topologically
Bergman. \pff

In order to see that the result is not void, we can exhibit an
action of $\Z$ by Lipschitz isomorphisms of $\R$ such that
$\{\ov{\Z\cdot x}\}_{x\in \R}$ does not partition $\R$. Namely,
let $T(x)=2 x$, whence $T^n(x)=2^nx$ for all $n\in \Z$. Then $T$
is a simple dilation of $\R$ and $0\in \ov{\Z\cdot x}$ for any
$x\in \R$.

We thus see from statement \ref{equicontinuity} that a Bergman group acting by H\"older maps actually acts equicontinuously. One might
wonder if this also holds if the group acts by, e.g., uniform homeomorphisms, but this is false. For example, $\ku S_\infty$,
which is the Bergman group par excellence, acts continuously on $\ca$, and thus by uniform homeomorphisms, but the action fails to
be equicontinuous. Moreover, there is no decomposition of $\ca$ into closed minimal pieces.

\section{Closure properties}
\begin{prop}\label{three space} Let $G$ be Polish
and $H\normal G$ a closed subgroup. It both $G/H$ and $H$ have
property (OB), so does $G$. Conversely, if $G$ has property (OB)
and $\phi:G\til K$ is a continuous homomorphism into a Polish
group $K$ with dense image, then $K$ has property (OB).
\end{prop}

\pf Let $V\subseteq G$ be an open neighbourhood of the identity in
$G$. Then $U=H\cap V$ is non-empty open in $H$, whence for some
finite set $A\subseteq H$ and $n\in \N$, $(UA)^n=H$. Therefore,
$(VA)^n\supseteq H$. Since the quotient mapping $\pi:G\til G/H$ is
continuous and open, $\pi[V]\subseteq G/H$ is open, non-empty, so
for some finite set $B\subseteq G$ and $m\in \N$,
$(\pi[V]\pi[B])^m=G/H$. Thus
\begin{equation}\begin{split}
\a g\in G\;\e v_1,\ldots,v_m\in V\;\e b_1,\ldots,b_n\in
B\;\;\big(\pi(g)=\pi(v_1b_1\cdots v_mb_m)\big)
\end{split}\end{equation}
and
\begin{equation}\begin{split}
&\a g\in G\;\e v_1,\ldots,v_m\in V\;\e b_1,\ldots,b_n\in B\;\e
v'_1,\ldots,v'_n\\
&\e a_1,\ldots,a_n\in A\;\; g=v_1b_1\cdots v_mb_mv'_1a_1\cdots
v'_na_n
\end{split}\end{equation}
So $G=(VB)^m(VA)^n$, showing that $G$ has property (OB).

Now, if $G$ has property (OB) and $K$ acts continuously by
isometries on a metric space, then this induces an action by $G$.
Thus, every $G$-orbit is bounded and as $\phi(G)$ is dense in $K$,
every $K$-orbit is bounded.\pff

\begin{prop} Suppose $\{G_n\}_\N$ are Polish groups. Then
$G=\prod_\N G_n$ has property (OB) if and only if each $G_n$ has
property (OB).\end{prop}

\pf Suppose each $G_n$ has property (OB) and assume $W_0\leq
W_1\leq\ldots\leq G$ is an exhaustive sequence of subsets with the
Baire property. As in the proof of Proposition \ref{cofinality} we
can suppose that each $W_n$ is a symmetric open neighbourhood of
the identity in $G$.

Thus there is $n\in \N$ such that
$$
\{1\}\times\ldots\times\{1\}\times \prod_{i>n}G_i\subseteq W_0
$$
and thus to prove that $G=W_n^k$ for some $n$ and $k$, it is
enough to prove that $G_0\times \ldots \times G_n$ is contained in
some $V_m^l$, where
$$
V_i=\{(g_0,\ldots,g_n)\in G_0\times \ldots \times G_n\del
(g_0,\ldots,g_n,1,1,\ldots)\in W_i\}
$$
But $V_0\subseteq V_1\subseteq \ldots\subseteq G_0\times \ldots
\times G_n$ is an increasing exhaustive sequence of open subsets,
and as $G_0\times \ldots \times G_n$ has property (OB) by
Proposition \ref{three space}, the result follows. The other
direction follows by Proposition \ref{three space}. \pff

In \cite{ber} Bergman poses the problem of whether (what we now
subsequently call) the Bergman property passes from a group to a
subgroup of finite index. In an announcement \cite{khe} A.
Kh\'elif states that this is indeed the case. We shall see that
the concept of induced representations also leads to this result
and, moreover, also solves the corresponding problem for Polish
groups.

\begin{prop}Let $G$ be a Polish group and $H\leq G$ a finite
index closed subgroup. The $G$ has property (OB) if and only if
$H$ has.\end{prop}

\pf First the easy direction. Assume $H$ has property (OB). Then
if $G$ acts continuously by isometries on some space $(X,d)$, so
does $H$ and this latter has bounded orbits. Letting
$g_1\,\ldots,g_n$ be representatives for the left cosets of $H$ in
$G$, we see that $G\cdot x=\for_ig_iH\cdot x$, which is a finite
union of bounded sets, and thus bounded.

For the other direction, consider first the abstract case of two
groups $G$ and $H$ with $H$ a finite index subgroup of $G$. Fix a
transversal $1\in T\subseteq G$ for the left cosets of $H$ in $G$.
Now assume that $H$ acts by isometries on a metric space $(X,d)$.
We define
$$
Y=\{ \xi:G\til X\del \a g\in G\;\a h\in H\; \xi(gh)=h\inv \xi(g)\}
$$
For example, if $x_0\in X$ is some fixed element, we can define
$\xi_0:G\til X$ by $\xi_0(ah)=h\inv x_0$ for all $h\in H$ and all
$a\in T$. Then clearly $\xi_0\in Y$. So $Y$ is non-empty.

We can now define the following metric $\partial$ on $Y$:
$\partial(\xi,\zeta)=\sup_{g\in G}d(\xi(g),\zeta(g))$. If we can
show that the supremum is finite, then this is clearly a metric.
But
\begin{equation}\begin{split}\label{finite sup}
\partial(\xi,\zeta)&=\sup_{g\in G}d(\xi(g),\zeta(g))\\
&=\sup_{a\in T, h\in H}d(\xi(ah),\zeta(ah))\\
&=\sup_{a\in T, h\in H}d(h\inv \xi(a),h\inv \zeta(a))\\
&=\sup_{a\in T, h\in H}d(\xi(a),\zeta(a))\\
&=\sup_{a\in T}d(\xi(a),\zeta(a))\\
&<\infty
\end{split}
\end{equation}
Where the last inequality holds since $T$ is finite. Now, let $G$
act on $Y$ by left translation
$$
(g\cdot \xi)(f)=\xi(g\inv f)
$$
This is an action by isometries.
\begin{equation}\begin{split}
\partial(g\cdot\xi,g\cdot\zeta)&=\sup_{f\in G}d((g\cdot\xi)(f),
(g\cdot\zeta)(f))\\
&=\sup_{f\in G}d(\xi(g\inv f),\zeta(g\inv f))\\
&=\sup_{f'\in G}d(\xi(f'),\zeta(f'))\\
&=\partial(\xi,\zeta)
\end{split}\end{equation}
Now, if $G$ has the Bergman property, there is a bounded orbit
$G\cdot \xi$ in $Y$. We now only need to see how this gives rise
to a bounded orbit for $H$ in $X$. So let $x_0=\xi(1)$ and notice
that for $h\in H$
\begin{equation}\begin{split}
d(x_0,h\cdot x_0)&=d(\xi(1),h\cdot\xi(1))\\
&=d(\xi(1),\xi(h\inv))\\
& =d(\xi(1),(h\cdot\xi)(1))\\
&\leq\sup_{g\in G}d(\xi(g),(h\cdot\xi)(g))\\
&=\partial(\xi,h\cdot\xi)\\
&\leq \textrm {diam}_\partial (G\cdot \xi)
\end{split}\end{equation}
This shows that the Bergman property passes to subgroups of finite
index.

For the case of Polish groups $G$ and $H$, with $H$ being a finite
index closed and thus clopen subgroup, we of course restrict our
attention to continuous $\xi$. Again we see that $\xi_0\in Y\neq
\tom$. We claim that the action of $G$ on $Y$ is separately
continuous. In the second variable this is trivial, as $G$ acts by
isometries. On the other hand, if we fix some $\xi\in Y$ and
suppose that $g_n\til g$ in $G$, then by Equation \ref{finite sup}
\begin{equation}\begin{split}
\partial(g_n\cdot\xi,g\cdot\xi)&
=\sup_{a\in T}d((g_n\cdot\xi)(a),(g\cdot\xi)(a)) =\sup_{a\in
T}d(\xi(g_n\inv a),\xi(g\inv a))\Lim{n\til\infty}0
\end{split}
\end{equation}
Thus, by Lemma \ref{joint}, the action of $G$ on $Y$ is continuous
and we can finish the proof as in the discrete case.\pff

These calculations will also help us prove
\begin{thm}\label{solvable}A solvable Bergman group is
finite.\end{thm}

This result reinforces our feeling that Bergman groups are not
something that can be constructed from below using only simple
construction methods, but are rather groups that have to be given
in one single step.

That infinite abelian groups might not be Bergman was suggested to
me by Pandelis Dodos. Actually, we shall see that there is a very
specific reason for this to be true.

\begin{defi}Let $H\leq K$ be abelian groups. $H$ is said to be
{\em pure} in $K$ if for every $n\geq 0$, $nK\cap H=nH$, i.e., if
some $x\in H$ is divisible by $n$ in $K$ then $x$ is divisible by
$n$ in $H$.

If, moreover, $K$ is a torsion group, we say that $H$ is a basic
subgroup of $K$ if (i) $H$ is pure in $K$, (ii) $H$ is a direct
sum of cyclic subgroups and (iii) $K/H$ is divisible.\end{defi}

A theorem of Kulikov (confer Robinson (4.3.4.) \cite{rob}) states
that any abelian torsion group contains a basic subgroup.

\begin{lemme}An infinite abelian group has a countably infinite
quotient group.
\end{lemme}

\pf Let $G$ be an infinite abelian group and $T$ its torsion
subgroup. Then $G/T$ is torsion free and hence embeds into a
direct sum of copies of $\Q$, $G/T\sqsubseteq  \oplus_{i\in I}\Q$.
If $G$ is not itself a torsion group, then the projection of $G/T$
onto one of these summands must be non-trivial and thus there is
some non-trivial quotient of $G/T$, and hence of $G$, which is
isomorphic to a subgroup of $\Q$. Clearly, this quotient is
countably infinite.

So we can assume that $G=T$, i.e., that $G$ is a torsion group.
Therefore, by the theorem of Kulikov, let $H$ be a basic subgroup
of $G$. Assume first that $|G/H|>1$. Then as $G/H$ is divisible,
by the structure theorem of divisible abelian groups, it is the
direct sum of isomorphic copies of quasicyclic groups. Since each
of these is countably infinite, by projecting $G/H$ surjectively
onto one of these we thus produce a countably infinite quotient of
$G$. Assume now that $G=H$. Then $G$ is a direct sum of cyclic
subgroups and, as $G$ is a torsion group, each of the summands is
finite. So let $G=\oplus_{i\in J} F_i$, where $F_i\neq \{0\}$ is a
finite group. As $G$ is infinite, so is $J$. Picking a countably
infinite subset $J_0\subseteq J$ and projecting $G$ onto
$\oplus_{i\in J_0}F_i$, we again end up with the countably
infinite quotient.\pff

Since no countably infinite group is Bergman, no infinite abelian
group is Bergman.

\

\noindent{\em Proof of Theorem \ref{solvable}.} Assume that $G$ is
solvable and that $1=G_0\normal G_1\normal\ldots\normal G_n=G$ is
an abelian series for $G$, i.e., that each quotient $G_{i+1}/G_i$
is abelian. Suppose that $G=G_n$ is Bergman and assume by
induction that $G_{i+1}$ is Bergman. Then also $G_{i+1}/G_i$ is
Bergman and, being abelian, it is also finite. Thus $G_i$ is a
finite index subgroup of a Bergman group and hence Bergman. So
this shows inductively that each of $G_i$ is Bergman, and so all
the quotients $G_{i+1}/G_i$ are finite. Therefore also $G$ is
finite.\pff

\section{Circle groups}

We shall first consider the homeomorphism group of the unit circle
$S^1$ and its model-theoretic counterpart, the automorphism group
of the countable dense circular order, ${\rm Aut}(\go C)$.

Let first $\pi : \R\til \R/\Z=S^1$ and let $d$ be the metric on
$S^1$ induced by the metric on $\R$. I.e., $d(x,y)=dist(\pi\inv
(x),\pi\inv (y))$. So $d$ takes values in $[0,1/2]$.

Let $\go C\subseteq S^1$ be a countable dense set, for
concreteness we can take $\go C=\pi[\Q]$, and ${\rm Aut}(\go C)$
the set of all permutations of $\go C$ that preserve the relation
$B\subseteq \go C^3$ defined as follows:

For $x,y,z\in \go C$ let $B(x,y,z)$ if and only if
\begin{itemize}
\item $x,y$ and $z$ are distinct.
\item Going clockwise along the unit circle $S^1$ from $x$ to $z$
one passes through $y$.
\end{itemize}
In this case, we say that $y$ is {\em between} $x$ and $z$.

\begin{nota}For $x_1,\ldots, x_n\in \go C$ write
$\circlearrowright x_1x_2\ldots x_n$ if for all $i<j<l$,
$B(x_i,x_j,x_l)$.\end{nota}

\begin{lemme}\label{produkt}
Suppose $x_1,\ldots,x_n,y_1,\ldots,y_n\in \go C$,
$\circlearrowright x_1\ldots x_ny_1\ldots y_n$. Then
$$
{\rm Aut}(\go C)={\rm Aut}(\go C,x_1,\ldots,x_n)\cdot {\rm
Aut}(\go C,y_1,\ldots,y_m)\cdot{\rm Aut}(\go C,x_1,\ldots,x_n)
$$
A similar statement holds for ${\rm Hom}_+(S^1)$, which is the
group of orientation preserving homeomorphisms of $S^1$.
\end{lemme}

\pf Let $g\in{\rm Aut}(\go C)$ be given and let $I\subseteq
\{1,\ldots,n\}$ be the set of $i$ such that $B(y_n,g(x_i),y_1)$.
Notice that $I$ is an interval, $I=\{i_0,i_0+1,\ldots, i_1\}$.

Pick some $f\in{\rm Aut}(\go C,y_1,\ldots,y_m)$ such that
$f(x_i)=g(x_i)$ for $i\in I$,
$$
B(y_n,f(x_i),x_1),\;\;\;\; B(y_n,f(x_i),f(x_{i_0}))
$$
for $i=1,\ldots, i_0-1$, and
$$
B(x_n,f(x_i),y_1),\;\;\;\; B(f(x_{i_1}),f(x_i),y_1)
$$
for $i=i_1+1,\ldots,n$.

Pick now
$$
h\in{\rm Aut}\big(\go
C,x_1,\ldots,x_n,f(x_{i_0}),\ldots,f(x_{i_1})\big)
$$
such that $hf(x_i)=g(x_i)$ for $i\notin I$. Then
$$
hf\begr\{x_1,\ldots,x_n\}=g\begr\{x_1,\ldots,x_n\}
$$
whence $(hf)\inv g\in{\rm Aut}(\go C,x_1,\ldots,x_n)$ and
$$
g=hf\cdot(hf)\inv g\in{\rm Aut}(\go C,x_1,\ldots,x_n)\cdot {\rm
Aut}(\go C,y_1,\ldots,y_m)\cdot{\rm Aut}(\go C,x_1,\ldots,x_n)
$$
The statement and the proof for ${\rm Hom}_+(S^1)$ is similar.\pff

\begin{defi}If $X$ is a set and ${\rm Sym}(X)$ the group of all
permutations of $X$, the {\em permutation group topology} on ${\rm
Sym}(X)$ has as open neighbourhood basis at the identity the sets
${\rm Sym}(X,x_1,\ldots,x_n)=\{g\in {\rm Sym}(X)\del
g(x_1)=x_1,\ldots, g(x_n)=x_n\}$, where $x_1,\ldots, x_n$ are any
elements of $X$.\end{defi}

Recall that on ${\rm Hom}_+(S^1)$ the two topologies of pointwise
and uniform convergence and the compact-open topology coincide.
Moreover, ${\rm Hom}_+(S^1)$ is a Polish group in this topology.

As ${\rm Aut}(\go C)$ is a closed subgroup of ${\rm Sym}(\N)=\ku
S_\infty$ in the permutation group topology, ${\rm Aut}(\go C)$ is
also Polish.

\begin{thm}${\rm Aut}(\go C)$ and ${\rm Hom}(S^1)$ are Bergman.
\end{thm}

\pf Let us first notice that for any $x\in \go C$, the groups
${\rm Aut}(\go C,x)$ and ${\rm Aut}(\Q,<)$ are naturally isomorphic. Thus if
$$
W_0\subseteq W_1\subseteq\ldots \subseteq {\rm Aut}(\go C)=\for_n
W_n
$$
then there are $n$ and $k$ such that ${\rm Aut}(\go C,x)\subseteq
W_n^k$. This follows from the result of Droste and G\"obel \cite{drogob} that ${\rm Aut}(\Q,<)$ is Bergman. Taking $g\in{\rm
Aut}(\go C)$ such that $x\neq g(x)$, we see by Lemma \ref{produkt} that
\begin{equation}\begin{split}
{\rm Aut}(\go C)&={\rm Aut}(\go C,x)\cdot {\rm Aut}(\go
C,g(x))\cdot{\rm Aut}(\go
C,x)\\
&={\rm Aut}(\go C,x)\cdot g{\rm Aut}(\go
C,x)g\inv \cdot{\rm Aut}(\go C,x)\\
&= W_n^kgW_n^kg\inv W_n^k\\
&= W_m^{5k}
\end{split}\end{equation}
for some sufficiently large $m$.

The same argument applies to ${\rm Hom}_+(S^1)$, using that for any $x\in S^1$ the groups
$$
{\rm Hom}_+(S^1,x)=\{g\in {\rm Hom}_+(S^1)\del g(x)=x\}
$$
and ${\rm Hom}_+(\R)$ are isomorphic. Again ${\rm Hom}_+(\R)$ is Bergman by the results of Droste and G\"obel, so by Lemma
\ref{produkt}, also ${\rm Hom}_+(S^1)$ is Bergman. Now the Bergman property for ${\rm Hom}(S^1)$ follows, as ${\rm Hom}_+(S^1)$ is a
subgroup of index $2$ in ${\rm Hom}(S^1)$.
\pff

\section{Groups of isometries}

\begin{defi} Let $(X,d)$ be a metric space and $G\leq {\rm
Iso}(X,d)$. We let $d_\infty$ denote the supremum metric on the
spaces $X^m$ induced by $d$. The group $G$ is said to be {\em
approximately oligomorphic} if for any $n\geq 1$ and $\eps>0$
there is a finite set $A\subseteq X^n$ such that $G\cdot A$ is
$\eps$-dense in $X^n$ with respect to $d_\infty$.\end{defi}

\begin{thm}\label{isometry groups} Let $(X,d)$ be a Polish
metric space and $G$ a closed subgroup of ${\rm Iso}(X,d)$ with
the topology of pointwise convergence. If $G$ is approximately
oligomorphic, then $G$ has property (OB).\end{thm}

\pf We need to show that $G$ is finitely generated of bounded
width over any non-empty open set $V\subseteq G$. So find $\ov
x=(x_1,\ldots,x_n)\in X^n$ and $\eps>0$ such that
$$
U=\{g\in G\del \a i\leq n\; d(x_i,gx_i)<\eps\}\subseteq VV\inv
$$
We claim that there is a finite set $B\subseteq X^n$ such that
$U\cdot B$ is $\frac{\eps}{2}$-dense in $X^n$. To see this, let
$A\subseteq X^n\times X^n$ be a finite set such that $G\cdot A$ is
$\frac{\eps}{2}$-dense in $X^n\times X^n$. Define $A'\subseteq A$
to be the set of $\ov a=(\ov a_1,\ov a_2)\in A$ such that for some
$g_{\ov a}\in G$, $d_\infty(\ov x, g_{\ov a}\ov
a_1)<\frac{\eps}{2}$. Finally, put $B=\{g_{\ov a}\ov a_2\del \ov
a\in A'\}$.

Then, if $\ov c\in X^n$, there are $\ov a=(\ov a_1,\ov a_2)\in A$
and $g\in G$ such that
$$
d_\infty\big((\ov x, \ov c),(g\ov a_1,g\ov
a_2)\big)<\frac{\eps}{2}
$$
In particular, $\ov a\in A'$ and thus
\begin{equation}\begin{split}
d_\infty(gg_{\ov a}\inv\ov x, \ov x)&= d_\infty(g_{\ov a}\inv\ov
x, g\inv\ov x)\\
&\leq d_\infty(g_{\ov a}\inv\ov x, \ov a_1)+d_\infty(\ov a_1,g\inv
\ov x)\\
&= d_\infty(\ov x, g_{\ov a} \ov a_1)+d_\infty(g\ov a_1, \ov x)\\
&<\eps \\
\end{split}\end{equation}
So $gg_{\ov a}\inv\in U$ and
$$
d_\infty(\ov c, gg\inv_{\ov a}\cdot g_{\ov a}\ov a_2)=d_\infty(\ov
c,g\ov a_2)<\frac{\eps}{2}
$$
Thus, $U\cdot B$ is $\frac{\eps}{2}$-dense in $X^n$. Let now
$B'\subseteq B$ be the set of $\ov b\in B$ such that
$$
d_\infty(U\cdot \ov b, G\cdot \ov x)<\frac{\eps}{2}
$$
So for some $h_{\ov b}\in G$ and $g_{\ov b}\in U$,
$$
d_\infty(g_{\ov b}\ov b, h_{\ov b}\ov x)<\frac{\eps}{2}
$$
Then, if $f\in G$, we can find $\ov b\in B$ and $g\in U$ such that
$d_\infty(f\ov x,g\ov b)<\frac{\eps}{2}$. Thus, $\ov b\in B'$ and
\begin{equation}\begin{split}
d_\infty(\ov x , f\inv gg_{\ov b}\inv h_{\ov b}\ov x)&=
d_\infty(g\inv
f\ov x,g_{\ov b}\inv h_{\ov b}\ov x)\\
&\leq d_\infty(g\inv f\ov x,\ov b)+
d_\infty(\ov b,g_{\ov b}\inv h_{\ov b}\ov x)\\
&=d_\infty(f\ov x,g\ov b)+d_\infty(g_{\ov b}\ov b, h_{\ov b}\ov
x)\\
&<\eps
\end{split}\end{equation}
So $f\inv gg_{\ov b}\inv h_{\ov b}\in U$ and if $H=\{h_{\ov b}\del
\ov b\in B'\}$, we see that $G=UU\inv HU\inv$.\pff

\begin{cor} Let $G$ be an oligomorphic closed subgroup of $\ku
S_\infty$. Then $G$ has property (OB).\end{cor}

\pf Notice that if we let $\N$ have the metric in which all points
have distance $1$, then $G$ is approximately oligomorphic as a
group of isometries exactly when it is oligomorphic.\pff

The {\em Urysohn metric space} $\U$ is the unique separable
complete metric space containing each finite metric space and such
that any isometry between finite subsets extends to a full
isometry of the space. This space, constructed by Urysohn
\cite{ury} is also characterised by being separable, complete and
satisfying the following extension property:

\noindent $(*)$: If $\phi:X\til \U$ is an isometric embedding of a
finite metric space $X$ into $\U$ and $Y=X\cap \{y\}$ is a one
point metric extension of $X$, then $\phi$ extends to an isometric
embedding of $Y$.

In the same manner, there is a Urysohn metric space of diameter
$1$, designated by $\U_1$,  which is the unique complete separable
metric space whose diameter is at most $1$ and satisfying the
extension property $(*)$, when $Y$ varies over metric spaces of
diameter at most $1$.

Similarly, one can construct variants of the Urysohn metric space,
where the metric takes values only in $\Q\cap [0,1]$. Thus, the
{\em rational Urysohn metric space of diameter $1$}, denoted by
$\Omega$, is the unique countable metric space whose metric takes
values in $\Q\cap [0,1]$ and satisfying the extension property for
$Y$, whose metric also takes values in $\Q\cap [0,1]$.

\begin{thm}\label{urysohn} Let $\U_1$ be the Urysohn metric space
of diameter $1$. ${\rm Iso}(\U_1)$ is approximately oligomorphic
and hence has property (OB).\end{thm}

For this proof we need some notions of metric theory. Let $D_n$ be
the set of $n\times n$ matrices $[a_{ij}]$ with entries in $[0,1]$
such that $d(i,j)=a_{ij}$ defines a pre-metric on
$\{1,\ldots,n\}$. Consider $D_n$ as a subset of $[0,1]^{n^2}$ with
the supremum metric $d_\infty$. Clearly, the triangle inequality
is a closed condition, so $D_n$ is compact.

We define also the following distance $d_1$ on $D_n$ \`a la Gromov
and Hausdorff (see Chapter 3, Gromov \cite{gro}):
\begin{displaymath}
d_1(A,B)=\min \big(trace(E)\del \left[\begin{array}{cc}
A & E\\
E^t & B\end{array}\right]\in D_{2n}\big)
\end{displaymath}
Notice that the infimum is indeed attained, as we are minimising
over a compact space. So if $A, B\in D_n$ are representing
pre-metrics $a$ and $b$ on $\{1,\ldots,n\}$ and $\{1',\ldots,
n'\}$ respectively (thus of diameter at most $1$), $d_1$ is the
minimum of $\sum_i c(i,i')$, where $c$ varies over all pre-metrics
on $\{1,\ldots,n,1',\ldots,n'\}$ of diameter at most $1$ agreeing
with $a$ on $\{1,\ldots,n\}$ and with $b$ on $\{1',\ldots,n'\}$.
Therefore, $d_1$ measures how far the spaces have to be from each
other, when they are both embedded into a metric space of diameter
at most $1$.

\begin{lemme}$2d_1\leq nd_\infty\leq nd_1$.
\end{lemme}

\pf Let $A, B\in D_n$ and let $a$ and $b$ be the corresponding
pre-metrics on $\{1,\ldots,n\}$. Assume that
$$
\delta=d_\infty(A,B)=\sup_{i,j}|a(i,j)-b(i,j)|
$$
and let $c$ be defined on $\{1,\ldots,n,1',\ldots,n'\}$ by
\begin{displaymath}
\begin{split}
c(i,j)=&a(i,j)\\
c(i',j')=&b(i,j)\\
c(i,j')=c(j',i)=& \min_l\big( a(i,l)+\delta/2+b(l,j)\big)
\end{split}\end{displaymath}
We claim that $c$ is a pre-metric and that $i\mapsto i$ and
$i\mapsto i'$ are isometric embeddings of the spaces given by $a$
and $b$ respectively. Clearly, the triangle inequality is
satisfied separately on $\{1,\ldots,n\}$ and on $\{1',\ldots,
n'\}$, and
\begin{equation}\begin{split}
c(i,k')+c(k',j)
&=\min_l\big(a(i,l)+\delta/2+b(l,k)\big)
+\min_p\big(b(k,p)+\delta/2+a(p,j)\big) \\
&=\delta+\min_l\big(a(i,l)+b(l,k)\big)
+\min_p\big(b(k,p)+a(p,j)\big) \\
&\geq\delta+\min_{l,p}\big(a(i,l)+b(l,p)+a(p,j)\big)\\
&\geq\min_{l,p}\big(a(i,l)+a(l,p)+a(p,j)\big)\\
&\geq a(i,j)\\
&=c(i,j)
\end{split}\end{equation}
Similarly, $c(i',j')\leq c(i',l)+c(l,j')$. And
\begin{equation}
\begin{split}
c(i,j')&\leq \min_k\big( a(i,k)+\delta/2+b(k,j)\big)\\
&\leq a(i,l)+\min_k\big( a(l,k)+\delta/2+b(k,j)\big)\\
&=c(i,l)+c(l,j')\\
\end{split}\end{equation}
Similarly, $c(i,j')\leq c(i,l')+c(l',j')$, so the triangle
inequality is verified. Unfortunately, $c$ does not necessarily
have diameter bounded by $1$, but this can be remedied by letting
$c'(x,y)=\min\{c(x,y),1\}$. Clearly, this does not affect the
distances on $\{1,\ldots,n\}$ and $\{1',\ldots, n'\}$ separately,
and only decreases other distances. So
$c'(i,i')=\frac{\delta}{2}$. Let now
\begin{displaymath}
C'=\left[\begin{array}{cc}
A & E\\
E^t & B\end{array}\right]\in D_{2n}
\end{displaymath}
be the matrix corresponding to $c'$, and notice that $d_1(A,B)\leq
trace(E)= n\frac{\delta}{2}=\frac{n}{2}d_\infty(A,B)$. Thus,
$d_1\leq\frac{n}{2}d_\infty$.

On the other hand, if the two sets $\{1,\ldots,n\}$ and
$\{1',\ldots, n'\}$ are very close to each other, pointwise, in
some common metric space, then the distance between $i$ and $j$
cannot differ very much from the distance between $i'$ and $j'$.
And in fact, $d_\infty\leq d_1$.\pff

\

\noindent{\em Proof of Theorem \ref{urysohn}:} Fix some $n\geq 1$
and $\eps>0$ and let $\partial$ be the metric on $\U_1$. As $D_n$
is compact, we can find some finite $\ku A\subseteq D_n$, which is
$\eps$-dense in the metric $d_1$. By the universality property of
the Urysohn metric space $\U_1$, this means that for any $\ov
x=(x_1,\ldots, x_n)\in \U_1^n$, there is $\ov y=(y_1,\ldots,
y_n)\in \U_1^n$ with distance matrix $A\in \ku A$, such that
$\partial_\infty(\ov x, \ov y)\leq \sum_i\partial(x_i,y_i)\leq
\eps$. So pick for each $A\in \ku A$ some $\ov z\in \U_1^n$ with
distance matrix $A$, and let $\AAA$ be the set of these. Then, if
$\ov x=(x_1,\ldots, x_n)\in \U_1^n$ there is $\ov y=(y_1,\ldots,
y_n)$ as above, and hence some $\ov z=(z_1,\ldots,z_n)\in \AAA$
isometric to $\ov y$. But then as $\U_1$ is ultrahomogeneous, we
see that $\ov y$ and $\ov z$ are in the same orbit of ${\rm
Iso}(\U_1)$, showing that ${\rm Iso}(\U_1)$ is approximately
oligomorphic.\pff

We will now show that if we consider the isometry group of the
rational Urysohn metric space of diameter $1$, $\Omega$, then we
actually get the Bergman property outright. The results here were
finally clear after a late night discussion at Caltech with Stevo
Todor\v cevi\'c.

Two tuples $\ov x$ and $\ov y$ in $\Omega$ are said to be {\em
uniformly of distance 1 from each other} if $d(x_i,y_j)=1$ for all
$i,j$.

\begin{lemme}\label{consistent extension} If $\ov x$ and $\ov y$
in $\Omega$ are uniformly of distance $1$ and some $\ov z$ in
$\Omega$ is given, then there are $\ov x'$ and $\ov z'$ such that
$(\ov x, \ov z, \ov y)$ and $(\ov x', \ov z', \ov y)$ are
isometric and $\ov x$ is uniformly of distance $1$ from both $\ov
x'$ and $\ov z'$.\end{lemme}

\pf Notice that the distances between $\ov x, \ov y, \ov z', \ov
x'$ are completely specified by the lemma, so we need only specify
the distances between $\ov z$ and $(\ov x',\ov z')$. We let $\ov
z$ be uniformly of distance $1$ from $\ov x'$ and put
$$
d(z_i,z'_j)=\min\{1,\inf_{y_l} d(z_i,y_l)+d(y_l,z_j)\}
$$
The triangle inequality holds, which  can be checked by hand, so
let us just give a few representative cases.

- Clearly, $d(x_i',y_j)\leq d(x_i',v)+d(v,y_j)$ for all $v\in \ov
x, \ov z, \ov y$, since $\ov x'$ is uniformly of distance $1$ from
all of $\ov x, \ov z, \ov y$. Moreover, it also holds for $v\in
\ov x',\ov z'$, since $(\ov x', \ov z', \ov y)$ is isometric to
$(\ov x, \ov z, \ov y)$ and thus is a metric space.

- Clearly, for all $v\in \ov x, \ov y,\ov x',\ov z'$,
$d(x_i,z'_j)\leq d(x_i,v)+d(v,z'_j)$, since in this case one of
the distances on the right hand side must equal $1$. So for
$v=z_j$ we have
\begin{equation}\begin{split}
d(x_i,z_k)+d(z_k,z'_j)&\geq\min\{1,\inf_{y_l}
d(x_i,z_k)+d(z_k,y_l)+d(y_l,z_j)\}\\
&\geq\min\{1,\inf_{y_l}
d(x_i,y_l)+d(y_l,z_j)\}\\
&\geq 1\\
&=d(x_i,z_j')
\end{split}\end{equation}

- Clearly, for all $v\in \ov x, \ov x',\ov y,\ov z$,
$d(z_i,y_j)\leq d(z_i,v)+d(v,y_j)$, since in the first two cases
one of the distances on the right hand side must equal $1$ and in
the last two cases it reduces to the triangle inequality on $(\ov
z, \ov y)$. And for $v=z'_j$ we have
\begin{equation}\begin{split}
d(z_i,z'_k)+d(z'_k,y_j)&
\geq\min\{1,\inf_{y_l}d(z_i,y_l)+d(y_l,z_k)+d(z_k,y_j)\}\\
&\geq\min\{1,d(z_i,y_j)\}\\
&\geq d(z_i,y_j)
\end{split}\end{equation}
\pff

\begin{lemme}\label{property (OB) for Omega}
Assume that $\ov x$ is a tuple in $\Omega$. Then there is some
$l\in {\rm Iso}(\Omega)$ such that
$$
{\rm Iso}(\Omega)=\big(l\cdot {\rm Iso}(\Omega,\ov x)\big)^4
$$
\end{lemme}

\pf Find some $\ov y$, isometric to $\ov x$ and uniformly of
distance $1$ from it. Let $l(\ov x)=\ov y$ and $l(\ov y)=\ov x$.
Then $l\cdot{\rm Iso}(\Omega,\ov x)\cdot l={\rm Iso}(\Omega,\ov
y)$. Let $g\in{\rm Iso}(\Omega)$ be any element and put $\ov
z=g(\ov x)$.

By Lemma \ref{consistent extension}, we can find $\ov x', \ov z'$
such that $(\ov x, \ov z, \ov y)$ and $(\ov x', \ov z', \ov y)$
are isometric and $\ov x$ is uniformly of distance $1$ from both
$\ov x'$ and $\ov z'$. Thus there is some $h\in{\rm
Iso}(\Omega,\ov y)$ such that $h(\ov x)=\ov x'$ and $h(\ov z)=\ov
z'$. Now, since $(\ov x, \ov z')$ and $(\ov x, \ov x')$ are
isometric, there is some $f\in{\rm Iso}(\Omega,\ov x)$ such that
$f(\ov z')=\ov x'$. And finally, as $(\ov y,\ov x')$ and $(\ov
y,\ov x)$ are isometric, we can find $k\in {\rm Iso}(\Omega,\ov
y)$ such that $k(\ov x')=\ov x$.

Therefore, $kfhg(\ov x)=kfh(\ov z)=kf(\ov z')=k(\ov x')=\ov x$ and
$kfhg\in{\rm Iso}(\Omega,\ov x)$. So
$$
g=h\inv f\inv k\inv (kfhg)\in\big({\rm Iso}(\Omega,\ov y)\cdot{\rm
Iso}(\Omega,\ov x)\big)^2=\big(l\cdot {\rm Iso}(\Omega,\ov x)\cdot
l\cdot {\rm Iso}(\Omega,\ov x)\big)^2.
$$
\pff

\begin{thm}\label{Bergman for Omega}
The isometry group of the rational Urysohn metric space of
diameter $1$, ${\rm Iso}(\Omega)$, is Bergman.\end{thm}

\pf The proof relies on the deep result of S. Solecki \cite{sol},
also independently announced by A.M. Vershik, that for any finite
rational metric space $X$ there is another finite rational metric
space $Y$ containing $X$ and such that any partial isometry of $X$
extends to a full isometry of $Y$.

First of all, we notice that this also implies the corresponding
result for rational metric spaces of bounded diameter $1$. For if
$X$ is of bounded diameter $1$, then we find first some $Y'$ (not
necessarily of bounded diameter $1$) extending $X$such that every
partial isometry of $X$ extends to a full isometry of $Y'$. Now,
if $d'$ is the metric on $Y'$, let $d$ be the metric given by
$d(y_0,y_1)=\min \{1, d'(y_0,y_1)\}$ and let $Y$ be the metric
space obtained. Then we see that the distances between points in
$X$ are preserved and thus $X$ is still a subspace of $Y$, and if
$f$ is an isometry of $Y'$ it is also an isometry of $Y$. Thus,
every partial isometry of $X$ extends to a full isometry of $Y$.

We now need the following concept, which will also be used in a
later section.

\begin{defi}Suppose $G$ is a Polish group and consider for each
finite $m\geq 1$ the diagonal conjugacy action of $G$ on $G^m$
given by
$$
g\cdot(h_1,\ldots,h_m)=(gh_1g\inv,\ldots,gh_mg\inv)
$$
$G$ is said to have {\em ample generics} if for each $m\geq 1$
there is a comeagre orbit in $G^m$ for this action.
\end{defi}

Notice that since $\Omega$ is countable, ${\rm Iso}(\Omega)$ is a
Polish group in the permutation group topology. In section 5 (A)
of Kechris and Rosendal \cite{kecros} it is shown how the above
extension property for finite rational metric spaces of bounded
diameter $1$ implies that ${\rm Iso}(\Omega)$ has ample generics
and the results of section 5 (F) of \cite{kecros} implies that a
Polish group with ample generics has the Bergman property if and
only if it has property (OB). Thus it is enough to show that ${\rm
Iso}(\Omega)$ is finitely generated of bounded width over any
non-empty open subset. But this follows from Lemma \ref{property
(OB) for Omega}.\pff

\section{A dense Bergman subgroup of the unitary group}
In the following $\ell_2$ will be the complex Hilbert space on the
countable orthonormal basis $(e_i)_{i\in\N}$ and with usual inner
product
$$
\langle \sum a_ie_i\del\sum b_ie_i\rangle=\sum a_i\ov{b_i}
$$
We will also fix a countable algebraically closed field
$\Q\subseteq\q\subseteq \C$ closed under complex conjugation. In
fact, it will only be essential that $\q$ is closed under square
root, and we could therefore work in some subfield of $\R$ too.
This would give similar results for the orthogonal group of the
real separable Hilbert space, but we shall be content with the
above setting. Notice first that $\q$ is dense in $\C$.

We let $\V$ be the $\q$-vector space with basis $(e_i)$ and notice
that $\V$ is a dense subset of $\ell_2$. The inner product
restricts to an inner product on $\V$ taking values in $\q$, as
$\q$ is a field closed under complex conjugation. Since $\q$ is
algebraically closed, the norm of an element of $\V$ also belongs
to $\q$. This will give us enough space to perform the usual tasks
of Gram-Schmidt orthonormalisation etc.

Our first result is the following extension property.

\begin{lemme}Let $T$ be a $\q$-linear isometry of $\V$. Then $T$
extends to a unique unitary operator on $\ell_2$.\end{lemme}

\pf Since $\V$ is dense in $\ell_2$ and $\ell_2$ is complete, any
isometry of $\V$ extends to a unique isometry of $\ell_2$. Hence
$T$ extends to an isometry of $\ell_2$ preserving the origin. A
simple argument shows that the extension is $\C$-linear.\pff

So the group $\ku U(\V)$ of $\q$-linear isometries of $\V$ can be
seen as a subgroup of $\ku U(\ell_2)$. It will be useful represent
unitary operators as infinite matrices with respect to the
canonical basis $(e_i)_{i\in\N}$. Since we are only considering
finite $\q$-linear combinations, this means that any row and any
column is eventually zero. The following operators in $\ku U(\V)$
are of particular interest.

\begin{defi} An operator $T\in \ku U(\V)$ is {\em finitary} if
it is the identity outside of a finite-dimensional subspace of
$\V$.\end{defi}

So the finitary operators are those that are supported on a
finite-dimensional subspace and hence can be represented as
\begin{displaymath}
\left[\begin{array}{cc}
A & 0\\
0 & I\end{array}\right]
\end{displaymath}
where $I$ is the infinite identity matrix and $A$ some finite
unitary matrix.

Clearly, the finitary operators form a subgroup of $\ku U(\V)$.

The unitary group $\ku U(\ell_2)$ naturally comes with the {\em
strong topology}, which makes it a Polish group. The strong
topology is the weakest topology that makes all the maps
$$
T\mapsto T(x)
$$
continuous, where $x$ varies over the elements of $\ell_2$.

Similarly, as $\V$ is countable, $\ku U(\V)$ is naturally
isomorphic to a subgroup if the group ${\rm Sym}(\V)$ of all
permutations of $\V$, with the Polish permutation group topology.
Moreover, $\ku U(\V)$ is easily seen to be closed in ${\rm
Sym}(\V)$ and hence is a Polish group itself. Notice that this
topology is stronger than the topology induced by $\ku U(\ell_2)$.

We need that the Gram-Schmidt orthonormalisation procedure can be
done in $\V$.

\begin{lemme}\label{gram}If $\W$ is any subspace of $\V$ and
$w_1,\ldots,w_n$ is an orthonormal set of vectors in $\W$, then
there is an orthonormal basis of $\W$ extending
$\{w_1,\ldots,w_n\}$.\end{lemme}

\pf Let $\{x_1, x_2,x_3,\ldots\}$ be a $\q$-vector space basis of
$\W$ such that $x_1=w_1,\ldots, x_n=w_n$. Now, define inductively
$y_m, w_m$ by
$$
y_{m+1}=x_m-\sum_{i=1}^m\langle x_{m+1}\del w_i\rangle w_i
$$
and notice that as $\langle x_{m+1}\del w_i\rangle\in \q$ also
$y_{m+1}\in \V$. Now, put
$$
w_{m+1}=\frac{y_{m+1}}{\|y_{m+1}\|}
$$
and again as $\|y_{m+1}\|\in \q$ (as $\q$ is algebraically
closed), $w_{m+1}\in \V$. So as usual, $\{w_1,w_2,\ldots\}$ is an
orthonormal basis of $\W$.\pff

\begin{lemme}\label{hrushovski}Suppose $S$ is a linear isometry
between finite-dimensional spaces $\W_0$ and $\W_1$. Then $S$
extends to a finitary operator $\tilde S$ in $\ku U(\V)$.
\end{lemme}

\pf This is clear from Lemma \ref{gram}. For choose an orthonormal
basis $v_1,\ldots, v_n$  for $\W_0$ and find some sufficiently big
$i$ such that $\W_0,\W_1\subseteq [e_1,\ldots,e_i]$. Then we can
extend $v_1, \ldots,v_n$ and $S(v_1),\ldots, S(v_n)$ respectively
to orthonormal bases $u_1,\ldots, u_i$ and $w_1,\ldots, w_i$  of
$[e_1,\ldots, e_i]$. Letting $\tilde S(u_j)=w_j$ for $j\leq i$ and
$\tilde S(e_j)=e_j$ for $j>i$, we have the result.\pff

In particular, if $\{v_1,\ldots,v_n\}$ and $\{u_1,\ldots,u_n\}$
are  orthonormal sets in $\V$, then there is a finitary operator
$F$ sending the ordered basis $\{v_1,\ldots,v_n\}$ to the ordered
basis $\{u_1,\ldots,u_n\}$.

We recall the following fact (see, e.g., Proposition 2.2. in
Kechris and Rosendal \cite{kecros})

\begin{prop}\label{comeagre} Let $G$ be a Polish acting continuously
on a Polish space $X$. Then the following are equivalent for any
$x\in X$:\\
(i) The orbit $G\cdot x$ is non-meager.\\
(ii) For every open neighbourhood $V\subset G$ of the identity,
$V\cdot x$ is somewhere dense.
\end{prop}

\begin{prop}\label{ample}$\ku U(\V)$ has ample generics.\end{prop}

\pf By abstract methods (see Truss \cite{truss} and Kechris and
Rosendal \cite{kecros}) it is enough to show that certain
amalgamation properties are satisfied, but in our case an outright
description of the comeagre orbits is not much longer, so we give
this instead.

So fix an $m\geq 1$. We need to construct $K_1, \ldots,K_m\in \ku
U(\V)$ such that the conjugacy class of the $m$-tuple
$(K_1,\ldots,K_m)$ is comeagre in $\ku U(\V)$.

By Lemma \ref{hrushovski} we see that for any $P_1,\ldots, P_m\in
\ku U(\V)$ and $v_1, \ldots, v_k\in \V$ there are finitary
operators $H_1, \ldots, H_m\in\ku U(\V)$ such that for all $t\leq
m, s\leq k$, $P_t(v_s)=H_t(v_s)$.

So list all $m$-tuples $\K=(K_1,\ldots,K_m)$ of unitary operators
on some common finite-dimensional subspace $\V_i=[e_1,\ldots,
e_i]$ as
$$
\K_1=(K^1_1,\ldots,K^1_m), \K_2=(K^2_1,\ldots,K^2_m),\ldots
$$
and let $a_n$ be the dimension of the space $[e_1,\ldots, e_i]$ on
which the $K_1^n, \ldots, K_m^n$ act. Let also
$b_n=\sum_{j=1}^na_j$. We suppose furthermore that each $\K_i$ is
repeated infinitely often.

We can now paste these operators together as
\begin{displaymath}M_t=
\left[\begin{array}{ccccc}
K_t^1 & &&\\
& K_t^2 &&\\
&& K_t^3&\\
&&&\ddots
\end{array}\right]\end{displaymath}
In other words, $M_1, \ldots, M_m$ are disjoint sums of unitary
operators on finite-dimensional spaces of the form
$[e_l,\ldots,e_k]$ such that each conjugacy type of $m$-tuples
appears infinitely often.

To see that the conjugacy type of $(M_1, \ldots, M_m)$ is comeagre
in $\ku U(\V)^m$, we show first that it is dense and non-meagre.
Thus, by Proposition \ref{comeagre}, it is enough to show that it
is dense and that for every $l\in \N$ the set
$$
\ku A=\big\{(T\inv M_1T, \ldots, T\inv M_mT)\del
T=\left[\begin{array}{cc}I_l&0\\
0&A\end{array}\right] \textrm{ for some } A \big\}
$$
is somewhere dense in $\ku U(\V)^m$, where $I_l$ is the $l\times
l$ identity matrix. To see this latter, find first some $b_i\geq
l$. We claim that $\ku A$ is dense in the open set of
$(P_1,\ldots,P_m)$ such that for every $t\leq m$,
\begin{displaymath}
P_t=\left[
\begin{array}{cccc}
K_t^1 & & \\
& \ddots  &  \\
& & K_t^i&  \\
& & & A_t
\end{array}\right]
=\left[
\begin{array}{cc}
M_t\begr _{\V_{b_i}} & \\
& A_t\\
\end{array}\right]
\end{displaymath}
for some $A_t$. For if $(P_1,\ldots,P_m)$ is above, then  the
tuple can be approximated arbitrarily well by a tuple of finitary
operators. So we can suppose that $P_1,\ldots,P_m$ are finitary
themselves. Assume that we want to approximate $P_1,\ldots, P_m$
on $\V_k=[e_1,\ldots, e_k]$, where $k>b_i$ is such that
$P_t(e_p)=e_p$ for all $t\leq m$ and $p\geq k$. Then we can find
$j>i$ such that $k=b_i+a_j$ and
\begin{displaymath}
P_t=\left[
\begin{array}{ccccc}
M_t\begr _{\V_{b_i}} & & \\
& K_t^j &\\
& & I
\end{array}\right]
\end{displaymath}
Find a unitary operator $T$ such that $T\begr[e_1,\ldots,e_{b_i}]=
I_{b_i}$ and $T$ sends the ordered basis
$\{e_{b_{j-1}+1},\ldots,e_{b_j}\}$ to
$\{e_{b_i+1},\ldots,e_{b_i+a_j}\}$. Then
\begin{displaymath}
TM_tT\inv= T\left[
\begin{array}{cccc}
K_t^1 & &\\
& K_t^2 &\\
&&\ddots
\end{array}\right]T\inv=
\left[
\begin{array}{ccccc}
M_t\begr _{\V_{b_i}} & & \\
& K_t^j &\\
& & B_t
\end{array}\right]
\end{displaymath}
for some $B_t$. Thus $TM_tT\inv$ agrees with $P_t$ on
$\V_k=[e_1,\ldots, e_{b_i+a_j}]$ for every $t=1,\ldots,m$.

A similar argument shows that the conjugacy class of $(K_1,\ldots,
K_m)$ is dense. But in fact, this also follows from the next
proposition. Thus, as there is a dense orbit, the diagonal
conjugacy action of $\ku U(\V)$ on $\ku U(\V)^m$ is generically
ergodic, i.e. any invariant Borel set is either meagre or
comeagre. So, as the conjugacy class of  $(K_1,\ldots, K_m)$ is
non-meagre, it must be comeagre.\pff

\begin{defi} A Polish group $G$ is said to have a {\em cyclically dense
conjugacy class} if there are elements $g,h\in G$ such that
$\{g^nhg^{-n}\}_{n\in \Z}$ is dense in $G$. We say that $G$ has an
{\em ample cyclically dense conjugacy class} if there is a $g\in
G$ and some infinite sequence $(h_k)_k\in G^\N$ such that the set
$\{(g^nh_kg^{-n})_k\}_{n\in \Z}$ is dense in $G^\N$.\end{defi}

\begin{prop}\label{cyclic} $\ku U(\V)$ has an ample cyclically
dense conjugacy class.\end{prop}

\pf Notice first that if  $G$ is a Polish group and for some $g\in
G$ there is some $(h_1,\ldots, h_m)\in G^m$ such that the set
$$\{(g^nh_1g^{-n},\ldots,g^nh_mg^{-n})\}_{n\in \Z}
$$
is dense in $G^m$, then set of such $(h_1,\ldots, h_m)$ is
certainly dense in $G^m$. Moreover, since it is also $G_\delta$,
it is comeagre. Thus, if for each $m\in \N$ there is such
$(h_1,\ldots, h_m)\in G^m$, then there is an infinite sequence
$(h_k)_k\in G^\N$ such that $\{(g^nh_kg^{-n})_k\}_{n\in \Z}$ is
dense in $G^\N$.

Therefore, we only need to find some unitary operator $S$ that
fills the r\^ole of $g$. For this we will consider instead a
biinfinite orthonormal basis $(e_i)_{i\in \Z}$ of $\V$ and let $S$
be the bilateral shift on this basis. Fix also some dimension $m$.

We can now take $H_t=I\oplus M_t$, where we let the identity $I$
act on $\V_-=[\ldots,e_{-2},e_{-1},e_0]$ and let $M_t$ be as in
the proof of Proposition \ref{ample} defined on
$\V_+=[e_1,e_2,e_3,\ldots]$. One easily sees that
$(S^{-n}H_1S^{n},\ldots, S^{-n}H_mS^{n})_{n\in \N}$ is dense in
$\ku U(\V)^m$. For suppose we wish to approximate some
$(P_1,\ldots,P_m)$, which we can suppose are finitary, on some
space $\W=[e_{-n},\ldots,e_n]$. Since each $P_t$ is finitary, we
can find $k>n$ such that each $P_t$ is on the form
\begin{displaymath}
P_t=\left[
\begin{array}{ccc}
I&&\\
&A_t & \\
&&I\end{array}\right]
\end{displaymath}
where $A_t$ is a $(2k+1)\times (2k+1)$ matrix acting on
$[e_{-k},\ldots,e_k]$. Now find some $j$ such that $K_t^j=A_t$ for
each $t\leq m$. Then we see that $H_t\begr [e_{b_{j-1}+1},\ldots,
e_{b_j}]$ and thus
$$
S^{-b_{j-1}-1-k}H_tS^{b_{j-1}+1+k}\begr[e_{-k},\ldots,e_k]=A_t
$$
Thus
$$
(S^{-b_{j-1}-1-k}H_1S^{b_{j-1}+1+k},\ldots,
S^{-b_{j-1}-1-k}H_mS^{b_{j-1}+1+k})
$$
agrees with $(P_1,\ldots,P_m)$ on $[e_{-k},\ldots,e_k]\supseteq
\W$.\pff

It follows from the results of Kechris and Rosendal (see
Proposition 5.16 \cite{kecros}) that a Polish group with ample
generics has property (OB) if and only if it is Bergman.

\begin{prop}$\ku U(\V)$ has the Bergman property.\end{prop}

\pf Since $\ku U(\V)$ has ample generics, it is enough to show
that it has property (OB). We show that $\ku U(\V)$ is finitely
generated of bounded width over any open neighbourhood $\mathsf U$
of the identity. So suppose $k$ is given such that $\mathsf U$
contains all operators fixing $\W_0=[e_1,\ldots,e_k]$ pointwise.
Find an operator $M$ such that
$M[\W_0]=\W_1=[e_{k+1},\ldots,e_{2k}]$ and notice that $M\mathsf U
M\inv$ contains all operators fixing $\W$ pointwise.  Let now
$T\in \ku U(\V)$ and find a finite dimensional space
$\HH_0\subseteq (\W_0\oplus \W_1)^\perp$ such that
$T[\W_0]\subseteq \W_0\oplus\W_1\oplus \HH_0$. Let $R_0\in\mathsf
U$ send  $\W_1$ into $(\W_0\oplus \W_1\oplus \HH_0)^\perp$  and
fix $\W_0\oplus\HH_0$ pointwise. Thus
\begin{displaymath}\begin{split}
R_0T[\W_0]\subseteq R_0[\W_0\oplus\W_1\oplus \HH_0] \subseteq
\W_0\oplus R_0[\W_1]\oplus \HH_0 \subseteq \W_1^\perp
\end{split}\end{displaymath}
We can therefore find some $R_1\in M\mathsf U M\inv$ such that
$R_1R_0T$ is the identity on $\W_0$ and $R_1$ fixes $\W_1$
pointwise. Hence, $R_1R_0T\in \mathsf U$ and $T\in R_0\inv R_1\inv
\mathsf U\subseteq \mathsf U\inv M\mathsf U\inv M\inv \mathsf U$.
Thus, $\ku U(\V)=\mathsf U\inv M\mathsf U\inv M\inv \mathsf
U$.\pff

\begin{lemme}$\ku U(\V)$ is dense in $\ku U(\ell_2)$.\end{lemme}

\pf It is enough to see that any unitary $T\in \ku U(\ell_2)$ can
be approximated arbitrarily well on any finite set of orthonormal
vectors. So suppose $x_1,\ldots, x_n$ is an orthonormal set and
$\eps>0$. By the continuity of the inner product, we can find
$\delta>$ such that if $v_1,\ldots,v_n,u_1,\ldots,u_n$ are
normalised vectors such that $\|x_i-v_i\|<\delta$ and
$\|T(x_i)-u_i\|<\delta$ for every $i$, then if $\hat
v_1,\ldots,\hat v_n$  and $\hat u_1,\ldots, \hat u_n$ are the
orthonormal bases obtained by applying the Gram-Schmidt
orthonormalisation process to $v_1,\ldots,v_n$ and
$u_1,\ldots,u_n$ respectively, we still have $\|x_i-\hat
v_i\|<\eps/2$ and $\|T(x_i)-\hat u_i\|<\eps/2$ for every $i$. Thus
choose $v_i, u_i\in \V$ as above and pick some $R\in \ku U(\V)$
sending the ordered basis $\hat v_1,\ldots,\hat v_n$ to the
ordered basis $\hat u_1,\ldots, \hat u_n$. Then
$$
\|T(x_i)-R(x_i)\|\leq \|T(x_i)-R(\hat v_i)\|+\|R(\hat
v_i)-R(x_i)\|\leq \eps
$$
and hence approximating $T$ on $x_1,\ldots,x_n$.\pff

So let us sum up the results so far.

\begin{thm}$\ku U(\V)$ has ample generics, an ample
cyclically dense conjugacy class and the Bergman property. Thus,
$\ku U(\ell_2)$ has property (OB) and an ample cyclically dense
conjugacy class.\end{thm}

Christopher Atkin \cite{atkin} has actually proved something a lot
stronger for the full unitary group, $\ku U(\ell_2)$, namely that
it has property (OB) in the norm topology.

We should mention that the existence of ample generics in a Polish
group has quite remarkable consequences for the structure of the
group, for example, it implies that any homomorphism from it into
a separable group is automatically continuous, and the group
cannot be covered by countably many non-open cosets (see Hodges,
Hodkinson, Lascar and Shelah \cite{hoholash}, Kechris and Rosendal
\cite{kecros}).

\begin{thm} Let $\ku U(\ell_2)$ act continuously by H\"older maps on
a complete metric space $(X,d)$. If the bilateral shift $S$
induces a relatively compact orbit on $X$, then $\ku U(\ell_2)$
fixes a point of $X$.\end{thm}

\pf First, we can evidently suppose that $X$ is in fact separable
and thus Polish. So $\ku U(\ell_2)$ acts continuously on $\ku
K(X)$. Moreover, if $g\in G$ is H\"older($\alpha$) with constant
$c$ on $(X,d)$, then $g$ is H\"older($\alpha$) with constant $c$
on $(\ku K(X),d_H)$, where $d_H$ is the Hausdorff metric. For
$$
d_H(K,L)=\max\big( \sup_{x\in K}d(x, L), \sup_{x\in L}d(K,x)\big)
$$
But
\begin{equation}\begin{split}
\sup_{x\in gK}d(x, gL)=\sup_{x\in K}\inf_{y\in L}d(gx,gy)\leq
\sup_{x\in K}\inf_{y\in L} c\cdot d(x,y)^\alpha=c\cdot (\sup_{x\in
K}\inf_{y\in L} d(x,y))^\alpha
\end{split}\end{equation}
and thus
\begin{equation}\begin{split}
d_H(gK,gL)&=\max\big( \sup_{x\in gK}d(x, gL), \sup_{y\in
gL}d(gK,y)\big)\\
&\leq \max\big( c\cdot (\sup_{x\in K}\inf_{y\in L} d(x,y))^\alpha,
c\cdot (\sup_{y\in L}\inf_{x\in K} d(x,y))^\alpha\big) \\
&=c\cdot \max \big(\sup_{x\in K}\inf_{y\in L} d(x,y),\sup_{y\in
L}\inf_{x\in K} d(x,y)\big)^\alpha\\
&=c\cdot d_H(gK,gL)^\alpha
\end{split}\end{equation}

Let now $\ku O=\ov {\{S^n\cdot x\}}_{n\in \Z}$ be compact and find
by the proof of Proposition \ref{cyclic} some $T\in \ku U(\ell_2)$
such that $\{S^nTS^{-n}\}_{n\in \Z}$ is dense in $\ku U(\ell_2)$.
Then
\begin{equation}\begin{split}
d_H(T\cdot \ku O,\ku O)&=d_H(TS^{-n} \ku O, S^{-n}\ku O)\\
&\leq c_nd_H(S^nTS^{-n}\ku O,S^nS^{-n}\ku O)^{\alpha_n}\\
&\leq c_nd_H(S^nTS^{-n}\ku O,\ku O)^{\alpha_n}
\end{split}\end{equation}
where $S^n$ is H\"older($\alpha_n$) with constant $c_n$. Using now
that $\ku U(\ell_2)$ has property (OB), we find a universal $N$
such that we can choose all the $\alpha_n\in [1/N,N]$ and $c_n\leq
N$. Picking a subsequence $n_i$ such that $S^{n_i}TS^{-n_i}\til
I$, we see that
$$
c_{n_i}d_H(S^{n_i}TS^{-n_i}\ku O,\ku O)^{\alpha_{n_i}}\til 0
$$
and thus $d_H(T\cdot \ku O,\ku O)=0$. Hence, $\ku O$ is both $S$
and $T$ invariant and thus also $\ku U(\ell_2)$-invariant.
Moreover, as $\ku O$ is compact and $\ku U(\ell_2)$ is extremely
amenable (this is a result of Gromov and Milman \cite{gromil}),
$\ku U(\ell_2)$ fixes a point of $\ku O$.\pff

\section{Groups of homeomorphisms}

\subsection*{Spheres}

\begin{thm}Let ${\rm Hom}(S^m)$ be the group of homeomorphisms of
the $m$-dimensional sphere with the topology of uniform
convergence. Then ${\rm Hom}(S^m)$ has property (OB).\end{thm}

\pf Let $d$ be the euclidean metric on $\R^{m+1}$ and $d_\infty$
the supremum metric on ${\rm Hom}(S^m)$, $d_\infty(g,f)=\sup_{x\in
S^m}d(gx,fx)$. We show that ${\rm Hom}(S^m)$ is finitely generated
of bounded width over any non-empty open subset $U$. So pick an
$\eps_0>0$ such that
$$
V=\{g\in {\rm Hom}(S^m)\del d_\infty(g,id)<3\eps_0\}\subseteq
UU\inv
$$
Let also $x_0=(1,0,0,\ldots,0)\in S^m$. Then for any
$\eps_0>\delta>0$ there is some homeomorphism $\phi_\delta$ of
$S^m$ such that $\phi_\delta(B_{\eps_0}(x_0))=B_\delta(x_0)$ and
$d_\infty(\phi_\delta,id)<\eps_0$. Moreover, there is an
involution homeomorphism $\psi$ of $S^m$ fixing $\partial
B_{\eps_0}(x_0)$ pointwise, while switching $int B_{\eps_0}(x_0)$
with $ext B_{\eps_0}(x_0)=S^m\setminus B_{\eps_0}(x_0)$. Finally,
let $\iota$ be the orientation inverting involution
$$
\iota(x_0,x_1,x_2,\ldots, x_m)=(x_0,-x_1,x_2,\ldots,x_m)
$$
We notice that $SO(m+1)$ is a compact subgroup of ${\rm
Hom}(S^m)$, so $\for_n (SO(m+1)\cap VV\inv)^n$ is an open subgroup
of $SO(m+1)$. But this latter is connected and compact, so
$SO(m+1)\subseteq (VV\inv)^k$ for some $k$.

\begin{claim}
${\rm Hom}(S^m)\subseteq (VV\inv)^k\{\iota, id\}V^2\psi V\psi
V\inv$
\end{claim}

To see this, let $g\in {\rm Hom}(S^m)$ and find $f\in SO(m+1)$
such that $fg(x_0)=x_0$. Then put
\begin{displaymath}
\hat f=\left\{ \begin{array}{c} f \;\;\;\textrm{ if } fg \textrm{
is orientation
preserving}\\
\iota f \;\;\;\textrm{ if } \iota fg \textrm{ is orientation
preserving}\end{array}\right.
\end{displaymath}
It follows that $\hat fg$ preserves the orientation and fixes
$x_0$. Therefore, by Lemma 3.1. of Glasner and Weiss
\cite{glawei}, which itself relies on the proof of the annulus
conjecture, there is some $\eps_0>\delta>0$ and a homeomorphism
$h$ of $S^m$ such that
\begin{equation}\begin{split}
&\a x\in S^m\;  (d(x,x_0)>\eps_0\til hx=\hat fgx)\\
&\a x\in S^m\;  (d(x,x_0)<\delta\til hx=x)
\end{split}\end{equation}
In particular,
$$
d_\infty(h,\hat fg)=\sup_{x\in S^m}d(hx,\hat fgx)=\sup_{x\in
B_{\eps_0}(x_0)}d(hx,\hat fgx)<3\eps_0
$$
So $\hat fgh\inv \in V$ and $g\in \hat f\inv Vh$. Thus,
$$
\phi\inv_\delta h\phi_\delta\begr B_{\eps_0}(x_0)=id
$$
and
$$
\psi\phi_\delta\inv h\phi_\delta\psi\begr ext B_{\eps_0}(x_0)=id
$$
In particular, $d_\infty(\psi\phi_\delta\inv
h\phi_\delta\psi,id)<3\eps_0$, so $\psi\phi_\delta\inv
h\phi_\delta\psi\in V$. Therefore, $h\in \phi_\delta\psi
V\psi\phi_\delta\inv$ and
\begin{equation}\begin{split}
g&\in \hat f\inv Vh\\
&\subseteq\hat f\inv V\phi_\delta\psi V\psi\phi_\delta\inv\\
&\subseteq (VV\inv )^k \cdot\{id, \iota\}\cdot V^2\psi V\psi V\inv
\end{split}\end{equation}\pff

\subsection*{The Hilbert cube}

Consider now the Hilbert cube $Q=[0,1]^\N$ and its homeomorphism
group ${\rm Hom}(Q)$ equipped with the topology of uniform
convergence. We let $d$ be the metric on $Q$ given by
$$
d((x_n),(y_n))=\sum_{n\in \N}\frac{|x_n-y_n|}{2^{n+1}}
$$
and $d_\infty$ the supremum metric on ${\rm Hom}(Q)$ given by
$$
d_\infty(f,g)=\sup_{\vec x\in Q}d(f(\vec x),g(\vec x)),
$$
which is right invariant.

\begin{thm}\label{hilbert cube}${\rm Hom}(Q)$, with the topology of
uniform convergence, has property (OB).\end{thm}

\pf Fix some open neighbourhood $V$ of the identity in ${\rm
Hom}(Q)$, which we can suppose is of the form
$$
V=\{g\in {\rm Hom}(Q)\del d_\infty(g,id)<2\eps\}
$$
for some $\eps>0$. Thus, if $n$ is sufficiently large such that
$2^{-n}<\eps$, then for any $f\in {\rm Hom}(Q)$ the does not
change the first $n$ coordinates of any $\vec x\in Q$, i.e.,
$$
f((x_0,x_1,\ldots,x_{n-1},x_n,x_{n+1},\ldots))=
(x_0,x_1,\ldots,x_{n-1},y_n,y_{n+1},\ldots)
$$
for all $\vec x\in Q$, we have $f\in V$.

\begin{claim}\label{transitivity}If $\frac{1}{k}<\eps$ , then
there is a finite set $\F\subseteq {\rm Hom}(Q)$ such that for
every $\vec x\in Q$, $\vec 0=(0,0,0,\ldots)\in \F V^{k+1}\cdot
\vec x$.\end{claim}

\noindent{\em Proof of claim:} Let $\frac{1}{2^n}<\eps$. For each
$s\in \{0,\frac{1}{2},1\}^n$ let $\vec
z_s=(s_0,s_1,\ldots,s_{n-1},0,0,\ldots)$. As $Q$ is homogeneous,
${\rm Hom}(Q)$ acts transitively on $Q$ (see van Mill, Theorem
6.1.6. \cite{vM}), and we can therefore find some $h_s\in{\rm
Hom}(Q)$ such that $h_s(\vec z_s)=\vec 0$ for each $s$. Let
$\F=\{h_s\del s\in \{0,\frac{1}{2},1\}^n\}$. So it is enough to
show that $\e s\; \vec z_s\in V^{k+1}\cdot \vec x$. So first use
the homogeneity of $Q$ to adjust the tail $(x_n,x_{n+1},\ldots)$
of $\vec x$ by some element of $V$ to become $(0,0,\ldots)$. This
can be done since a homeomorphism leaving the first $n$
coordinates invariant belongs to $V$. Now we can subsequently
adjust each of the first $n$ coordinates (leaving the tail
invariant) to be equal to either $0,\frac{1}{2}$ or $1$. For this
operation it is enough to use a product of at most $k$ elements of
$V$.\pff

Now, it follows from Brouwer's fixed point Theorem that any
homeomorphism of $Q$ fixes a point. Thus, up to a conjugate by an
element of the set $\F V^{k+1}$ from Claim \ref{transitivity}, any
homeomorphism of $Q$ fixes $\vec 0$. As we wish to show that ${\rm
Hom}(Q)$ is finitely generated of bounded width over $V$, we can
suppose that any homeomorphism fixes $\vec 0$.

\begin{claim}(Glasner and Weiss \cite{glawei})\label{fix-points}
If $f\in {\rm Hom}(Q)$ fixes $\vec 0$, then there is some
$\delta>0$ and $g\in{\rm Hom}(Q)$ such that $d_\infty(g,f)<\eps$
and $g\begr _{B_\delta(\vec 0)}=id$.\end{claim}

\noindent{\em Proof of claim:} Pick $\delta>0$ sufficiently small
such that $\sup_{\vec x\in B_\delta(\vec 0)} d(f(\vec x), \vec
x)<\eps$. As  both $\partial B_\delta(\vec 0)$ and
$\partial[f"B_\delta(\vec 0)]$ are $Z$-sets, we can extend the
homeomorphism $f\inv :\partial[f"B_\delta(\vec 0)]\til \partial
B_\delta(\vec 0)$ to a homeomorphism $h\in {\rm Hom}(Q)$
satisfying $d_\infty(h,id)<\eps$ (see van Mill, Theorem 6.4.6.
\cite{vM}). Thus, $d_\infty(hf,f)=d_\infty(h,id)<\eps$ and we can
let
\begin{displaymath}
g(\vec y)=\left \{
\begin{array}{cc} \vec y \;\;\;&\textrm {if}\; \vec y\in
B_\delta(\vec 0)\\
hf(\vec y)&\textrm{otherwise}\end{array}\right.
\end{displaymath}
\pff

\begin{claim}\label{compressing}
For any $g\in {\rm Hom}(Q)$ and $0<\delta<\eps$ such that $g\begr
_{B_\delta(\vec 0)}=id$, there is $h\in V^2$ such that $h\inv
gh\begr _{{[0,\eps[}^{n+1}\times [0,1]^\N}=id$.
\end{claim}

\noindent{\em Proof of claim:} Notice that
${{[0,\delta[}^{n+1+l}\times [0,1]^\N}\subseteq B_\delta(\vec 0)$
for some $l>0$. Moreover, it is not hard to see that
$[0,\eps[\times [0,1]^{l}$ is homeomorphic to ${[0,\delta[}^{l+1}$
by some function $a$, which is a homeomorphism of $[0,1]^{l+1}$. Thus,
$$
h_0=id_{[0,1]^n}\otimes a\otimes id_{[0,1]^\N}:Q\til Q
$$
belongs to $V$ and sends
$$
{[0,1]^n}\times \Big([0,\eps[\times [0,1]^{l}\Big)\times{[0,1]^\N}
$$
to
$$
{[0,1]^n}\times{[0,\delta[}^{l+1}\times{[0,1]^\N}.
$$
Now, let $h_1:Q\til Q$ be a homeomorphism that moves the set
${[0,\eps[}^n\times [0,1]^\N$ to ${[0,\delta[}^n\times [0,1]^\N$,
preserves all coordinates $\geq n$ and $d_\infty (h_1,id)<\eps$.
Then $h_0,h_1\in V$ and $h=h_1h_0$ moves ${{[0,\eps[}^{n+1}\times
[0,1]^\N}$ to ${{[0,\delta[}^{n+1+l}\times [0,1]^\N}$. \pff

Now, let $\iota\in {\rm Hom}([0,1])$ be an involution
homeomorphism that fixes $\eps$ and switches $0$ and $1$. Define
$i\in {\rm Hom}(Q)$ by
$$
i(x_0,x_1,\ldots,x_n,x_{n+1},\ldots)=(\iota(x_0),\iota(x_1),\ldots,
\iota(x_{n}),x_{n+1},x_{n+2},\ldots).
$$
Then $i$ interchanges ${{[0,\eps[}^{n+1}\times [0,1]^\N}$ and
${{]\eps,1]}^{n+1}\times [0,1]^\N}$.

We can now conclude our result. For suppose $f\in{\rm Hom}(Q)$.
Then, up to a conjugate by an element of $\F V^{k+1}$ we can
suppose that $f$ fixes $\vec 0$. By Claim \ref{fix-points}, we can
find $g\in Vf$ and $0<\delta<\eps$ such that
$g\begr_{B_\delta(\vec 0)}=id$. So pick by Claim \ref{compressing}
some $h\in V^2$ such that $h\inv gh\begr _{{[0,\eps[}^{n+1}\times
[0,1]^\N}=id$. But then $ih\inv ghi\begr_{{]\eps,1]}^{n+1}\times
[0,1]^\N}=id$ and hence $ih\inv ghi\in V$. All in all, this shows
that $f\in (\F V^{k+1})\inv V\inv V^2iViV^{-2}\F V^{k+1}$.\pff

Since any Polish group is a closed subgroup of ${\rm Hom}(Q)$ (see
Uspenski\u\i  \; \cite{usp} or the exposition in Kechris
\cite{kec}), we have

\begin{cor} Any Polish group is topologically isomorphic to a
closed subgroup of a Polish group with property (OB).\end{cor}

\section{Actions on trees}
\subsection*{Comeagre conjucagy classes}

We give first a simple proof of a result of Dugald Macpherson and
Simon Thomas. We will actually prove a result stronger than
theirs, for which we need some basic computations by M. Culler and
J.W. Morgan \cite{cm}. Note first that if a group $G$ acts by
isometries on an $\R$-tree $T$, then each $g\in G$ has associated
a characteristic non-empty subtree $T_g$ of $T$, which either is
the set of points fixed by $g$ (in which case $g$ is called {\em
elliptic}) or a line on which $g$ acts by translation (in which
case $g$ is called {\em hyperbolic}). We let $\|g\|=\inf(r\in
\R_+\del \e x\in T\; d(x,g\cdot x)=r)$. This infimum is in fact
attained as shown in \cite{cm}. Thus, $g$ is elliptic if and only
if $\|g\|=0$.

The interested reader should consult the very readable article by
Culler and Morgan for more information on the general theory of
group actions on $\R$-trees.

\begin{lemme}\label{cm1}\cite{cm} Suppose $g$ and $h$ are isometries
of an $\R$-tree $T$ . If $T_g\cap T_h$ is empty, then
$$
\|gh\|=\|g\|+\|h\|+2dist(T_g,T_h)
$$
\end{lemme}

\begin{lemme}\label{cm2}\cite{cm} Let $g$ and $h$ be hyperbolic
isometries of an $\R$-tree $T$ such that $T_g\cap T_h\neq \tom$.
Then
$$
\max(\|gh\|,\|gh\inv\|)=\|g\|+\|h\|
$$
\end{lemme}
From Lemma \ref{cm1} follows the following important special case.

\begin{thm}(Serre's Lemma) Suppose $g$, $h$ and $gh$ are elliptic
isometries of an $\R$-tree $T$. Then $T_g\cap T_h\neq
\tom$.\end{thm}

\begin{thm}\label{macpherson}(D. Macpherson and S. Thomas for combinatorial trees
\cite{mactho}.) Suppose $G$ is a Polish group with a comeagre
conjugacy class $C$ acting by isometries on an $\R$-tree $T$. Then
every element of $G$ is elliptic.\end{thm}

\pf We claim that $\|\cdot\|$ is constantly $0$ on $C$. Assume
towards a contradiction that this is not the case. Notice first
that $\|\cdot\|$ is conjugacy invariant, so constant on $C$. Pick
$g,h\in C$ such that also $gh,gh\inv \in C$. By Lemma \ref{cm1},
if $T_g\cap T_h=\tom$ then
$$
\|gh\|=\|g\|+\|h\|+2dist(T_g,T_h)>\|g\|
$$
contradicting that $\|\cdot\|$ is constant on $C$. So $T_g\cap
T_h\neq \tom$, whence by Lemma \ref{cm2},
$$
\max(\|gh\|,\|gh\inv\|)=\|g\|+\|h\|>\|g\|
$$
again contradicting that $\|\cdot\|$ is constant on $C$ and thus
proving the claim.

Assume now that $f$ is an arbitrary element of $G$ and pick
$g,h\in C$ such that $f=hg$. Then
$$
C\inv\cap C\cap hC\inv\cap gC\inv\cap fC\inv\neq \tom
$$
so we can find $k_0,k_1,k_2,k_3\in C$ with
$k_0=hk_1\inv=gk_2\inv=fk_3\inv$, $k_0\inv\in C$, i.e.,
$k_0k_1=h$, $k_0k_2=g$ and $k_0k_3=f=k_0k_1k_0k_2$.

Notice that $k_0,k_1,k_0k_1\in C$, $k_1,k_0k_2,k_1k_0k_2=k_3\in C$
and $k_0\inv, k_0k_2,k_0\inv k_0k_2=k_2\in C$, so applying Serre's
Lemma to each of these three situations, we have $T_{k_0}\cap
T_{k_1}\neq \tom$, $T_{k_1}\cap T_{k_0k_2}\neq \tom$ and
$T_{k_0}\cap T_{k_0k_2}\neq \tom$. The three trees $T_{k_0}$,
$T_{k_1}$ and $T_{k_0k_2}$ therefore intersect pairwise, and thus
there is some $x$ in their common intersection. But then clearly
$f\cdot x=k_0k_1k_0k_2\cdot x=x$, whence $f$ is elliptic. \pff

Notice that if a group $G$ acts by automorphisms on a tree $T$,
then the action extends to the tree $T'$ obtained from $T$ by
adding a midpoint on every edge. Moreover, the action on $T'$ is
{\em without inversion}, i.e., there are no vertices $a\neq b$ in
$T'$ such that $\{a,b\}$ is an edge and $g\cdot a=b$, $g\cdot b=a$
for some $g\in G$. We also see that $G$ fixes a vertex if $T'$ if
and only if $G$ fixes either a vertex or an edge of $T$. Thus, to
see that a group has property (FA) it is enough to show that any
action without inversion on a tree has a fixed vertex.

The proof of Theorem \ref{macpherson} translates word for word into the corresponding proof
for $\Lambda$-trees (a generalisation of $\R$-trees with a metric
taking values in an arbitrary ordered abelian group). The only
thing that has to be checked is that the appropriate lemmas are
true also in this setting. Well, here they are. In the following,
$\Lambda$ is a fixed ordered abelian group and $(X,d)$ a given
$\Lambda$-tree. We define the norm of elements of $G$ in the same
manner as for actions on $\R$-trees.

\begin{lemme}\label{chis1}(Lemma 2.1.11 in \cite{chis}) Suppose
$X_1,\ldots,X_n$ are subtrees of $X$ such that $X_i\cap X_j\neq
\tom$ for all $i,j$. Then $X_1\cap\ldots\cap X_n\neq
\tom$.\end{lemme}

\begin{lemme}\label{chis2}(Lemma 3.2.2 in \cite{chis})
Suppose $g$ and $h$ are isometries of $(X,d)$, which are not
inversions, such that $T_g\cap T_h=\tom$. Then
$$
\|gh\|=\|g\|+\|h\|+2dist(T_g,T_h)
$$
\end{lemme}

\begin{lemme}\label{chis3}(Lemmas 3.2.3 and 3.3.1 in \cite{chis})
Suppose $g$ and $h$ are hyperbolic isometries of $(X,d)$ such that
$T_g\cap T_h\neq \tom$. Then
$$
\max(\|gh\|,\|gh\inv\|)=\|g\|+\|h\|
$$
\end{lemme}
From Lemma \ref{chis2} we have again a version of Serre's Lemma.
\begin{lemme}Suppose $g$, $h$ and $gh$ are elliptic
isometries of a $\Lambda$-tree $(X,d)$. Then $T_g\cap T_h\neq
\tom$.\end{lemme}

\begin{thm} Suppose $G$ is a Polish group with a comeagre conjugacy
class $C$ acting by isometries and without inversion on a
$\Lambda$-tree $(X,d)$. Then every element of $G$ is
elliptic.\end{thm}

\subsection*{Dense conjugacy classes}

\begin{lemme}\label{continuity} Suppose a topological group
$G$ acts continuously and without inversion on a tree $T$, i.e.,
such that the stabilisers of vertices in $T$ are open in $G$. Then
$\norm{\cdot}: G\til \N$ is continuous.\end{lemme}

\pf Suppose first that $\norm{g}=0$. Then for some $a\in T$,
$g\cdot a=a$, i.e. $g\in G_a$ and $G_a$ is an open neighbourhood
of $g$ on which $\norm{\cdot}$ is constantly $0$.

Now, suppose $\norm{g}=n>0$. Then by a theorem of Tits
(Proposition 24 in \cite{serre}) there is a line $\ell_g=(a_i\del
i\in \Z)$ in $T$ such that $g\cdot a_i=a_{i+n}$ for all $i$. Now,
if $h\in G$ is elliptic, then for any $a\in T$, $h$ fixes the
midpoint of the geodesic from $a$ to $h\cdot a$. So if $h\cdot
a_0=g\cdot a_0=a_n$ then $n=2m$, $m>0$ and $h\cdot a_m=a_m\neq
g\cdot a_m$, by uniqueness of the geodesic. Hence if
$$
U=\{f\in G\del f\cdot a_0=g\cdot a_0,\ldots,f\cdot a_n=g\cdot
a_n\}
$$
then $U$ is an open neighbourhood of $g$ containing only
hyperbolic points of norm $\leq n$.

Moreover, if $h$ is hyperbolic of norm $k<n$, then $\ell_h$ would
contain exactly the $k+1$ midpoints of the arc $a_0, a_1,\ldots,
a_n$ from $a_0$ to $h\cdot a_0=a_n$. So for some $0<i<n$,
$$
dist_T(a_i,h\cdot a_i)=k\neq n
$$
which is a contradiction. So $U$
only contains hyperbolic points of norm $n$.\pff

\begin{prop}\label{dense}Suppose $G$ is a Polish group with a
dense conjugacy class, which is not the union of a countable
sequence of proper open subgroups. Then whenever $G$ acts
continuously and without inversion on a tree $T$, it fixes a
vertex of $T$. In other words, $G$ has property (topFA).\end{prop}

\pf Notice first that $\norm{\cdot}$ is conjugacy invariant and
continuous, so must be constantly $0$ on $G$. I.e. every element
of $G$ is elliptic. So if $G$ does not fix a vertex, it fixes an
end $\alpha=(a_0,a_1,\ldots)\subseteq T$ (Tits, Exercise 2, page
66 \cite {serre}). But then $G=\for_nG_{(a_n,a_{n+1},\ldots)}$,
where $G_{(a_n,a_{n+1},\ldots)}$ is the pointwise stabiliser of
the set $\{a_n,a_{n+1},\ldots\}$. Since these subgroups are
closed, almost all of them must be open, as $G$ satisfies Baire's
category theorem. And as $G$ is not the union of a countable chain
of proper open subgroups, $G=G_{(a_n,a_{n+1},\ldots)}$ for some
$N$, contradicting that $G$ did not fix a vertex.\pff

S. Solecki \cite{sol} has shown that the isometry group of the
rational Urysohn metric space, ${\rm Iso}(\U_\Q)$, with the
permutation group topology, has ample generics and a cyclically
dense conjugacy class. Moreover, in Kechris and Rosendal
\cite{kecros} it is shown that Polish groups with ample generics
and a cyclically dense conjugacy class cannot be written as the
union of a countable chain of proper subgroups. So this means that
${\rm Iso}(\U_\Q)$ has property (FA). Moreover, V.G. Pestov \cite{pestov} shows that ${\rm Iso}(\U)$ has no non-trivial continuous representations by isometries in a reflexive Banach space, so in particular it has property (FH). However, this does not solve the corresponding problem for ${\rm Iso}(\U_\Q)$.

{\em Address:}{ Christian Rosendal, Department of Mathematics, 273 Altgeld Hall MC-382,

1409 W. Green Street, Urbana, IL 61820, USA.}

{\em Email:} rosendal@math.uiuc.edu


\begin{thebibliography}{99}
\bibitem{atkin}Atkin, Christopher J. {\em Boundedness in uniform
spaces, topological groups, and homogeneous spaces.} Acta Math.
Hungar. 57 (1991), no. 3-4, 213-232.



\bibitem{becker}Becker, Howard {\em Polish group actions:
dichotomies and generalized elementary embeddings.} J. Amer. Math.
Soc. 11 (1998), no. 2, 397--449.


\bibitem{beckec}Becker, Howard \& Kechris, Alexander S. {\em The
descriptive set theory of Polish group actions.} London
Mathematical Society Lecture Note Series, 232. Cambridge
University Press, Cambridge, 1996. xii+136 pp.



\bibitem{bek} Bekka, Bachir {\em Kazhdan's property (T) for the
unitary group of a separable Hilbert space.} Geom. Funct. Anal. 13
(2003), no. 3, 509--520.



\bibitem{behava} Bekka, Bachir; de la Harpe, Pierre \& Valette, Alain
{\em Kazhdan's property (T)}, forthcoming book 2003.


\bibitem{ber}Bergman, George M. {\em Generating infinite symmetric
groups.} To appear in the "Bulletin of the London Mathematical
Society".




\bibitem{chis} Chiswell, Ian {\em Introduction to $\Lambda$-trees.}
World Scientific Publishing Co., Inc., River Edge, NJ, 2001.
xii+315 pp.



\bibitem{corn} de Cornulier, Yves {\em Uncountable groups with
property (FH).} To appear in "Communications in Algebra".



\bibitem{cm}Culler, Marc \& Morgan, John W. {\em Group actions
on $\R$-trees.} Proc. London Math. Soc. (3) 55 (1987), no. 3,
571--604.


\bibitem{drogob}Droste, Manfred \& G\"obel, R\"udiger {\em Uncountable
cofinalities of permutation groups.} To appear in the "Procedings
of the London Mathematical Society".


\bibitem{drohol}Droste, Manfred \& Holland, W. Charles {\em Generating
automorphism groups of chains.} To appear in "Forum Math.".

\bibitem{glawei}Glasner, Eli \& Weiss, Benjamin {\em The topological
Rohlin property and topological entropy.} Amer. J. Math. 123
(2001), no. 6, 1055--1070.


\bibitem{gro}Gromov, Misha {\em Metric structures for Riemannian and
non-Riemannian spaces.} Based on the 1981 French original. With
appendices by M. Katz, P. Pansu and S. Semmes. Translated from the
French by Sean Michael Bates. Progress in Mathematics, 152.
Birkhäuser Boston, Inc., Boston, MA, 1999. xx+585 pp.

\bibitem{gromil}Gromov, Misha \& Milman, Vitaly D. {\em A topological
application of the isoperimetric inequality.} Amer. J. Math. 105
(1983), no. 4, 843--854.



\bibitem{hej1} Hejcman, Jan {\em Boundedness in uniform spaces and
topological groups.} Czechoslovak Math. J. 9 (84) 1959 544--563.



\bibitem{hej2}Hejcman, Jan {\em On simple recognizing of bounded
sets.} Comment. Math. Univ. Carolin. 38 (1997), no. 1, 149--156.



\bibitem{greg} Hjorth, Greg {\em Classification and orbit equivalence
relations.} Mathematical Surveys and Monographs, 75. American
Mathematical Society, Providence, RI, 2000. xviii+195 pp.

\bibitem{hoholash}Hodges, Wilfrid; Hodkinson, Ian; Lascar,
Daniel \& Shelah, Saharon {\em The small index property for
$\omega$-stable $\omega$-categorical structures and for the random
graph.} J. London Math. Soc. (2) 48 (1993), no. 2, 204--218.




\bibitem{kec} Kechris, Alexander S. {\em Classical descriptive set
theory.} Graduate Texts in Mathematics, 156. Springer-Verlag, New
York, (1995). xviii+402 pp.


\bibitem{kecros} Kechris, Alexander S. \& Rosendal, Christian
{\em Turbulence, amalgamation and generic automorphisms of
homogeneous structures}, preprint (2004) (available at
www.math.caltech.edu/$\sim$rosendal).

\bibitem{khe}Kh\'elif, Anatole {\em On the Bergman property},
announcement, author's email address: khelif@logique.jussieu.fr.



\bibitem{koptits}Koppelberg, Sabine \& Tits, Jacques {\em Une
propri\'et\'e des produits directs infinis de groupes finis
isomorphes.} (French) C. R. Acad. Sci. Paris Sér. A 279 (1974),
583--585.



\bibitem{mactho} Macpherson, Dugald \& Thomas, Simon {\em Comeagre
conjugacy classes and free products with amalgamation.} preprint
2003.

\bibitem{vM}van Mill, Jan {\em Infinite-dimensional topology.
Prerequisites and introduction.} North-Holland Mathematical
Library, 43. North-Holland Publishing Co., Amsterdam, 1989.
xii+401 pp.

\bibitem{mil} Miller, Benjamin D. {\em Full groups, classification,
and equivalence relations.} Dissertation, UC Berkeley 2004,
(available at http://www.math.ucla.edu/~bdm/papers.html).


\bibitem{pestov}Pestov, V.G. {\em The isometry group of the Urysohn space as a L\'evy group.} preprint 2005.



\bibitem{rob}Robinson, Derek J. S. {\em A course in the theory
of groups.} Second edition. Graduate Texts in Mathematics, 80.
Springer-Verlag, New York, 1996. xviii+499 pp.



\bibitem{saxl}\textsc{Saxl, J.; Shelah, S.; Thomas, S.} {\em Infinite products of finite simple groups.}  Trans. Amer. Math. Soc.  348  (1996),  no. 11, 4611--4641.




\bibitem{serre}Serre, Jean-Pierre {\em Trees.} Translated from the
French original by John Stillwell. Corrected 2nd printing of the
1980 English translation. Springer Monographs in Mathematics.
Springer-Verlag, Berlin, 2003. x+142 pp.



\bibitem{shelah} Shelah, Saharon
{\em On a problem of Kurosh, J\'onsson groups, and applications.}
Word problems, II (Conf. on Decision Problems in Algebra, Oxford, 1976), pp. 373--394,
Stud. Logic Foundations Math., 95,
North-Holland, Amsterdam-New York, 1980. 




\bibitem{sol} Solecki, S\l awomir {\em Extending partial
isometries.} to appear in Israel Journal of Mathematics.

\bibitem{thomas} \textsc{Thomas, S.} {\em Infinite products of finite simple groups. II.}  J. Group Theory  2  (1999),  no. 4, 401--434.


\bibitem{tol1} Tolstykh, Vladimir {\em Infinite dimensional
general linear groups are groups of universally finite width},
preprint 2004.

\bibitem{tol2} Tolstykh, Vladimir {\em On Bergman's property for
the automorphism group of relatively free groups}, preprint 2004.

\bibitem{truss}Truss, John K. {\em Generic automorphisms of
homogeneous structures.} Proc. London Math. Soc. (3) 65 (1992),
no. 1, 121--141.

\bibitem{ury}Urysohn, P. {\em Sur un espace m\'etrique universel.}
Bull. Sci. Math. 51 (1927), 43-64, 74-90.

\bibitem{usp} Uspenski\u\i, Vladimir V. {\em A universal topological
group with a countable basis.} (Russian) Funktsional. Anal. i
Prilozhen. 20 (1986), no. 2, 86--87.

\end{thebibliography}
\end{document}